\documentclass[review,onefignum,onetabnum]{siamart190516}

\usepackage{amsfonts}
\usepackage{graphicx}
\usepackage{algorithmic}
\usepackage{arydshln}

\usepackage{mathtools}
\usepackage{multicol}
\usepackage{jlcode_JX}

\usepackage{subcaption}

\ifpdf
  \DeclareGraphicsExtensions{.eps,.pdf,.png,.jpg}
\else
  \DeclareGraphicsExtensions{.eps}
\fi

\newsiamremark{remark}{Remark}
\newsiamremark{hypothesis}{Hypothesis}
\crefname{hypothesis}{Hypothesis}{Hypotheses}
\newsiamthm{claim}{Claim}

\headers{MetaFEM: A Generic FEM Solver By Meta-expressions}{J. Xie, K. Ehmann, and J. Cao}

\title{MetaFEM: A Generic FEM Solver By Meta-expressions\thanks{
\funding{This work was funded by the Department of Defense Vannevar Bush Faculty Fellowship N00014-19-1-2642}}}

\author{Jiaxi Xie
\and Kornel Ehmann
\and Jian Cao\thanks{Department of Mechanical Engineering, Northwestern University, Evanston, IL (\email{jcao@northwestern.edu}).}}

\usepackage{amsopn}

\ifpdf
\hypersetup{
  pdftitle={MetaFEM: A Generic FEM Solver By Meta-expressions},
  pdfauthor={J. Xie, K. Ehmann, and J. Cao}
}
\fi

\begin{document}

\maketitle

\begin{abstract}
  Current multi-physics Finite Element Method (FEM) solvers are complex systems in terms of both their mathematical complexity and lines of code.
  This paper proposes a skeleton generic FEM solver, named MetaFEM, in total about 5,000 lines of Julia code, which translates 
  generic input Partial Differential Equation (PDE) weak forms 
  into corresponding GPU-accelerated simulations with a grammar similar to FEniCS or FreeFEM. 
  Two novel approaches differentiate MetaFEM from the common solvers:
  (1) the FEM kernel is based on an original theory/algorithm which explicitly processes meta-expressions, as the name suggests, 
  and (2) the symbolic engine is a rule-based Computer Algebra System (CAS), i.e., the equations are rewritten/derived 
  according to a set of rewriting rules instead of going through completely fixed routines, supporting easy customization by developers. 
  Example cases in thermal conduction, linear elasticity and incompressible flow are presented to demonstrate utility. 
\end{abstract}



\section{Introduction} 

Multi-physics FEM solvers are an important component of the fundamental infrastructure 
in a wide range of modern engineering and academic fields.
Most common solvers, if not all, are complex. Their complexity can be roughly expressed in term of Lines Of Code (LOC), i.e.,  
each solver typically has $>10^5$ LOC, as seen by the examples in \cref{tab:FEM_packages}. Naturally,
it is tempting to have a compact, flat skeleton software that encompasses a large portion of major functions, i.e., solves a wide range of PDE systems.
\begin{table}[tbhp]
  {\footnotesize
    \caption{Codebase size of selected open-source multi-physics FEM solvers, more FEM package names (but not data) can be simply found, 
    for example, in the list of Wikipedia \cite{WikiList}.} \label{tab:FEM_packages}
  \begin{center}
    \begin{tabular}{|c|c|c|} \hline
      Name & Module and Version & LOC (Comments Included)  \\ \hline
      FreeFEM \cite{FreeFEM} & 4.9 & 166805 lines of C++ in /src  \\ \hdashline
      GOMA \cite{GOMA} & 6.2 & 423285 lines of C in /src and /include \\ \hdashline
      Elmer \cite{Elmer} & 9.0 & 375613 lines of Fortran or C in /fem/src \\ \hdashline
      FEniCS, \cite{FEniCS_2012} & DOLFIN \cite{DOLFIN_2010} 0.3.1 & 49895 lines of C++ in /cpp/dolfinx, \\ 
      DOLFIN-toolchain & & 31004 lines of Python in /python/dolfinx \\ 
            & FFC \cite{FFC_2006} 0.3.1 & 10035 lines of Python in /ffcx \\
            & UFL \cite{UFL_2014} 2019.1.0 & 23255 lines of Python in /ufl \\ \hdashline
      MOOSE \cite{MOOSE} & 0.9.0 & 175745 lines of C in /framework/src, \\ 
            & & (/modules uncounted)\\
      \hline
    \end{tabular}
  \end{center}
  }
\end{table}

Different multi-physics solvers have very different architectures and levels of generality. 
This paper focuses on a small specific type, namely generic FEM solvers, 
which process completely generic PDE weak forms. 
That is, the solver merely uses tensor symbols with (free or dumb) indices added, multiplied, or algebraically operated, 
without higher level concepts such as Navier-Stokes equation(s), Dirichlet boundary condition(s), or even advection and diffusion. 
More directly, the script of the physics should look just like the mathematical expressions someone would write on a blackboard, 
for example, in a continuum mechanics class. 

To process generic PDE weak forms, the solver needs a symbolic processor and a FEM kernel.
The symbolic processor parses the input script and outputs a reorganized intermediate representation. 
The FEM kernel links the intermediate representation to the mesh and generates 
the well-known final linear system $\mathbf{K}\mathbf{x} = \mathbf{d}$.
For a clear illustration, the FEM kernel will be introduced before the symbolic processor.

A FEM kernel can have very different formulations, but its essential function is relatively fixed, i.e., the linear  
system $\mathbf{K}\mathbf{x} = \mathbf{d}$ is uniquely described by:
\begin{enumerate}
  \item PDE weak forms (domain, boundary, stabilization);
  \item Linearization, e.g., complete gradient of nonlinear terms;
  \item Element type and order;
  \item Sufficient quadrature order or the same numerical integration scheme;
  \item Temporal discretization scheme; and
  \item The sequence (numbering) of variables.
\end{enumerate}
which the kernel "simply" assembles. 

A symbolic processor, however, can be fundamentally improved by applying the theory of rewriting \cite{RW_system_1999}. 
A general symbolic operation is a mapping between two symbolic expressions 
and can be described by a third expression composed of the former two, denoted by a rewriting rule,
which essentially establishes a mapping from operations to data. 
A symbolic processor is no more than a collection of elemental symbolic operations, and can be 
practically generated by parsing a set of simple rules, which is just a list of data. 
Compared to the classical approaches in which the input expressions are parsed, automatically 
differentiated, re-organized and assembled by hard-coded routines, e.g., those in FEniCS and 
FreeFEM, the two most well-developed generic FEM solvers to the author’s best knowledge, 
the proposed rule-based formulation is flatter in structure, smaller in size and much more extensible for customization 
by simply adding more rules, resulting in a practical CAS with the above functions 
plus simplification on the expression level in $\approx$ 1,500 LOC.

With the above discussions, this paper proposes MetaFEM, a compact open-source generic FEM solver of $\approx$ 5,000 lines of Julia \cite{Julia_2017} code, 
with the only dependence on CUDA.jl \cite{CUDA_Julia_2019} for the GPU interface and some other Julia's intrinsic libraries. 
Compared to not only the above classical solvers in Fortran/C/C++/Python but also the relatively newer solvers in Julia like Gridap.jl \cite{Gridap},   
MetaFEM comes with an original architecture, i.e.: 
(1) the FEM kernel is based on a straightforward formulation by directly processing meta-expressions, and 
(2) uses a rule-based CAS by applied symbolic rewriting as the symbolic processor. 

The paper is organized as follows. The FEM kernel is described in
\cref{sec:FEM_kernel}, the rewriting system in \cref{sec:RW_system}, numerical
results are given in \cref{sec:results}, with the conclusions following in \cref{sec:conclusions}.

\section{The FEM kernel}
\label{sec:FEM_kernel}
To begin with, the notation for two different kinds of collections need to be clarified for the rest of the paper.
First, a physical vector is denoted by an arrow over the symbol, e.g., $\vec{u}$, and a component is
subscripted in English $i, j, ...$ like $u_i$, with the Einstein summation convention, e.g., $u_iu_i\coloneqq\sum_{i=1}^{dim}u_iu_i$. 
Second, a general, variable-sized collection is enclosed in curly brackets like $\{\phi_\alpha\}$, and a component is
subscripted in Greek $\alpha, \beta, ...$ like $\phi_\alpha$, without the Einstein summation convention.
The collection size, which is also the maximum index, is denoted by a hat, e.g., $\hat{\alpha}\coloneqq\alpha_{max}$.

\subsection{Theory}
For illustration simplicity, we begin with a PDE system with one variable $\phi$ on one single compact manifold $\Omega$ 
with boundary $\partial\Omega$. From the engineering perspective, $\Omega$ is simply a workpiece assigned with some known physics.

Then, each PDE in $\Omega$ is limited to the following meta-expression:
\begin{equation}
    \label{eq:meta}
    \mathcal{L}(\phi) = \mathcal{L}^a(\partial_t^{\nu_1}D_1\phi, ..., \partial_t^{\nu_{\lambda}}D_{\lambda}\phi, ..., \partial_t^{\nu_{\hat{\lambda}}}D_{\hat{\lambda}}\phi) = 0
\end{equation}
where $\mathcal{L}$ is the overall operator of $\phi$ with the algebraic operator $\mathcal{L}^a$, 
which has the arity $\hat{\lambda}$, i.e., the number of operands. Each operand can be addressed by its index $\lambda = 1,2,..., \hat{\lambda}$. 
The $\lambda^{th}$ operand has the ${(\nu_\lambda)}^{th}$ order  
temporal differential operator $\partial_t^{\nu_{\lambda}}$ and the spatial differential operator $D_{\lambda}$.

A weak solution $\phi^w$ (superscript $w$ for weak) is a function which satisfies:
\begin{equation}
    \int_{\Omega}\mathcal{L}(\phi^w)\delta\overline{\phi}^w = 0 
\end{equation}
under an arbitrary smooth test function $\overline{\phi}^w$. Each specific $\overline{\phi}^w$ is independent of $\phi^w$, but they will be always discussed in pairs,
so the overline is used here and later to distinguish the symbols while emphasizing their duality.

The FEM, or more generally the minimum weighted residual method(s), 
is to find a discretized solution $\phi^h$ (superscript $h$ following the customary notation in literature) which satisfies 
\begin{equation}
    \label{eq:single_wf}
    \int_{\Omega}\mathcal{L}(\phi^h)\delta\overline{\phi}^h = 0 
\end{equation}
under the arbitrary discretized test function $\overline{\phi}^h$. Each discretized function is a function which can be decomposed into 
the weighted sum of some predetermined interpolation functions:
\begin{align}
    \partial_tD\phi^h(\vec{x},t) & = \sum_{\alpha=1}^{\hat{\alpha}}(DN_\alpha(\vec{x}))(\partial_t\phi_{\alpha}(t)) \\  
    \delta(D\overline{\phi}^h(\vec{x},t)) & = \sum_{{\overline{\alpha}}=1}^{\hat{\alpha}}(D\overline{N}_{\overline{\alpha}}(\vec{x}))\delta\overline{\phi}_{{\overline{\alpha}}} 
\end{align}
where $N_\alpha$ is the $\alpha^{th}$ shape function and $\overline{N}_{\overline{\alpha}}$ is the $\overline{\alpha}^{th}$ test function, with 
$\alpha, \overline{\alpha} = 1,2, ..., \hat{\alpha}$. 
Shape functions and test functions are two predetermined collections of interpolation functions of size $\hat{\alpha}$,
which span the base space and the dual space, respectively.
The order of accuracy is determined by the complete polynomial order of the base space 
while the order of conservation is determined by the complete polynomial order of the dual space.
When the dual space contains the base space, 
the scheme yields optimal convergence. 
Therefore, the most natural design is to choose $N_\alpha = \overline{N}_\alpha$ for each $\alpha$, as in classical FEM, 
Discontinuous Galerkin (DG), etc., where the overline is usually simply omitted since no distinguishing is needed. 

In classical FEM, each $\alpha$ is (the index of) a specific mesh node and 
$\phi_\alpha$ is the value of the corresponding physical variable at that node. 
In general, each $N_\alpha$ represents a deformation mode 
with $\phi_{\alpha}$ representing the corresponding deformation amplitude, 
named by a control point value, but may not be equal to any physical variable at a geometric point.

By directly extending \cref{eq:meta}, 
a generic PDE system can be represented as a sum of bilinear forms, i.e., as:
\begin{equation}
    d = \overbrace{(\cdot,\cdot)_{\Omega}+...+(\cdot,\cdot)_{\Omega}}^{\text{domain physics}}+
    \overbrace{(\cdot,\cdot)_{\partial\Omega}+...+(\cdot,\cdot)_{\partial\Omega}}^{\text{boundary conditions}}+
    \overbrace{(\cdot,\cdot)_{\Omega}+...+(\cdot,\cdot)_{\Omega}}^{\text{numerical modifications}}
\end{equation}
where the overall weak form $d$ consists of the domain physics, the boundary conditions and the numerical modifications 
(in addition to the domain physics such as stabilizations), while each
\begin{equation}
    (D_0\overline{\phi}, \mathcal{L}^a(..., \partial_t^{\nu_{\lambda}}D_{\lambda}\phi, ...))_{\Omega'}\coloneqq\int_{\Omega'}\mathcal{L}^a(..., \partial_t^{\nu_{\lambda}}D_{\lambda}\phi, ...){\delta}D_0\overline{\phi}
\end{equation}
is a single bilinear-form on $\Omega' = \Omega \text{ or } \partial\Omega$ with the dual word $D_0\overline{\phi}$ and the base term
$\mathcal{L}^a(..., \partial_t^{\nu_{\lambda}}D_{\lambda}\phi, ...)$. A word is a symbol
optionally with subscripts and a term is an expression tree formed by words.
The details of words and terms will be discussed in \cref{sec:RW_system} and are skipped here.

Numerically, each bilinear form is explicitly approximated by numerical integration as:
\begin{align}
    (D_0\overline{\phi}, \mathcal{L}^a &(..., \partial_t^{\nu_{\lambda}}D_{\lambda}\phi, ...))_{\Omega} \approx \\ 
    \nonumber \sum_{\gamma}w_{\gamma}^{\text{itg}}\mathcal{L}^a &
    (...,\sum_{\alpha = 1}^{\hat{\alpha}}(D_\lambda N_{\alpha})|_{\vec{x}_{\gamma}^{\text{itg}}}(\partial_t^{\nu_\lambda}\phi_{\alpha}), ...)
    \sum_{{\overline{\alpha}} = 1}^{\hat{\alpha}}(D_0\overline{N}_{{\overline{\alpha}}})|_{\vec{x}_{\gamma}^{\text{itg}}}\delta\overline{\phi}_{{\overline{\alpha}}}
\end{align}
where $w_{\gamma}^{\text{itg}}, \vec{x}_{\gamma}^{\text{itg}} $ are the weight and position of the $\gamma^{th}$ numerical integration point respectively.

With \cref{eq:single_wf}, classical FEM minimizes each discretized component 
\begin{equation}                               
    d_{{\overline{\alpha}}}(\{\phi_{\alpha}\}) \coloneqq\frac{\partial{d}}{\partial(\delta\overline{\phi}_{{\overline{\alpha}}})} 
\end{equation}
by the Newton-Raphson method. For a static problem, with an initial (guess or given) set of control point values $\{\phi_{\alpha}^0\}$, 
for each sub-step $n$ (in the overall one single timestep) one has:
\begin{equation}
    \label{eq:Newton-Raphson}
    0=d_{{\overline{\alpha}}}|_{\{\phi_{\alpha}^{n+1}\}} \approx \sum_{\alpha=1}^{\hat{\alpha}}\frac{\partial{d_{{\overline{\alpha}}}}}{\partial\phi_{\alpha}^n}|_{\{\phi_{\alpha}^{n}\}}
    (\phi_{\alpha}^{n+1}-\phi_{\alpha}^{n})+d_{{\overline{\alpha}}}|_{\{\phi_{\alpha}^{n}\}}
\end{equation}
where we denote
\begin{equation}
    \label{eq:stiffness}
    K_{{\overline{\alpha}}\alpha}^n\coloneqq\frac{\partial{d_{{\overline{\alpha}}}}}{\partial\phi_{\alpha}^n}
\end{equation}
to be the $({\overline{\alpha}}, \alpha)^{th}$ component of the tangent stiffness of sub-step $n$, 
resulting in the linear system 
$\sum_{\alpha=1}^{\hat{\alpha}}K^n_{{\overline{\alpha}}\alpha}(\phi_{\alpha}^{n+1}-\phi_{\alpha}^{n}) = -d_{{\overline{\alpha}}}|_{\{\phi_{\alpha}^{n}\}}$,
which is exactly the well established $\mathbf{Kx} = \mathbf{d}$ in classical literature. \cref{eq:stiffness} will
be explicitly calculated in \cref{subsec:FEM_Algor}.

For a dynamic problem, the above process needs to be extended according to the temporal discretization, 
and we choose the generalized-$\alpha$ scheme \cite{generalized_alpha_1993} 
with the maximum temporal derivative order $\hat{\nu} = 2$ for a simple but practical example. 

At each timestep $m$, sub-step $n$, the collection of
\begin{equation}
    \{\phi_{\alpha,m}^{n+1}, ...\}\coloneqq\{\phi_{\alpha,m}^{n+1}, u_{\alpha,m}^{n+1}, a_{\alpha,m}^{n+1}\}
\end{equation}
needs to be determined as the basic variable, velocity and acceleration, where the term "basic" specifically emphasizes the symbol without time derivative,
with the following constraints:
\begin{equation}
    \Delta\phi_{\alpha, m}^n = \phi_{\alpha, m}^{n+1}-\phi_{\alpha, m}^n, \quad
    \Delta{u}_{\alpha, m}^n = u_{\alpha, m}^{n+1}-u_{\alpha, m}^n, \quad
    \Delta{a}_{\alpha, m}^n = a_{\alpha, m}^{n+1}-a_{\alpha, m}^n
\end{equation}
\begin{equation}
    \Delta\phi_{\alpha, m}^n = \Delta{t}(u_{\alpha, m}^0+b_1\Delta{u}_{\alpha, m}^n), \qquad
    \Delta{u}_{\alpha, m}^n = \Delta{t}(a_{\alpha, m}^0+b_2\Delta{a}_{\alpha, m}^n)
    \label{eq:time_constraints}
\end{equation}
\begin{equation}
    \widetilde{\phi}_{\alpha, m}^n = \phi_{\alpha, m}^{0}+c_1\Delta\phi_{\alpha, m}^n, \quad
    \widetilde{u}_{\alpha, m}^n = u_{\alpha, m}^{0}+c_2\Delta{u}_{\alpha, m}^n, \quad
    \widetilde{a}_{\alpha, m}^n = a_{\alpha, m}^{0}+c_3\Delta{a}_{\alpha, m}^n
    \label{eq:time_effective}
\end{equation}
where $\{\Delta\phi_{\alpha,m}^{n}, ...\}$ are the $(\hat{\nu} + 1) \times \hat{\alpha}$ incremental values 
linked by $\hat{\nu} \times \hat{\alpha}$ constraints defined by the predefined constraint parameters $b_1, b_2$, 
so that the only $\hat{\alpha}$ degrees of freedom are exclusively the basic variables $\{\phi_{\alpha}\}$ but not its time derivatives, 
as in \cref{eq:time_constraints}.
Meanwhile, $\{\widetilde{\phi}_{\alpha, m}^n, ...\}$ are the $(\hat{\nu} + 1) \times \hat{\alpha}$ effective values at which the residues are evaluated with 
$c_1, c_2 \text{ and } c_3$ being the relaxation parameters interpolating between the values at the last time step and the current (incremental) values, 
as in \cref{eq:time_effective}. 

By extending \cref{eq:Newton-Raphson}, the corresponding linear system becomes:
\begin{align}
    \nonumber 0={}&d_{{\overline{\alpha}}}|_{\{\widetilde{\phi}_{\alpha, m}^{n+1}, ...\}} \qquad\text{for } {\overline{\alpha}} = 1,2,...,\hat{\alpha}
    \\ \nonumber \approx{}&\sum_{\alpha=1}^{\hat{\alpha}}
    (\frac{\partial{d_{{\overline{\alpha}}}}}{\partial\widetilde{\phi}_{\alpha}^n}
    \frac{\partial\widetilde{\phi}_{\alpha}^n}{\partial\phi_{\alpha}^n}+
    \frac{\partial{d_{{\overline{\alpha}}}}}{\partial\widetilde{u}_{\alpha}^n}
    \frac{\partial\widetilde{u}_{\alpha}^n}{\partial\phi_{\alpha}^n}+ 
    \frac{\partial{d_{{\overline{\alpha}}}}}{\partial\widetilde{a}_{\alpha}^n}
    \frac{\partial\widetilde{a}_{\alpha}^n}{\partial\phi_{\alpha}^n})|_{\{\widetilde{\phi}_{\alpha, m}^{n}, ...\}}(\phi_{\alpha,m}^{n+1}-\phi_{\alpha,m}^{n})+d_{{\overline{\alpha}}}|_{\{\widetilde{\phi}_{\alpha, m}^n, ...\}}
    \\ ={}&\sum_{\alpha=1}^{\hat{\alpha}}
    (c_1\frac{\partial{d_{{\overline{\alpha}}}}}{\partial\widetilde{\phi}_{\alpha}^n} +
    \frac{c_2}{b_1\Delta{t}}\frac{\partial{d_{{\overline{\alpha}}}}}{\partial\widetilde{u}_{\alpha}^n} +
    \frac{c_3}{b_1b_2(\Delta{t})^2}\frac{\partial{d_{{\overline{\alpha}}}}}{\partial\widetilde{a}_{\alpha}^n})|_{\{\widetilde{\phi}_{\alpha, m}^{n}, ...\}}
    (\phi_{\alpha,m}^{n+1}-\phi_{\alpha,m}^{n})
    \label{eq:gen_alpha}
    +d_{{\overline{\alpha}}}|_{\{\widetilde{\phi}_{\alpha, m}^n, ...\}}
\end{align}
and when $d_{{\overline{\alpha}}}|_{\{\widetilde{\phi}_{\alpha, m}^n, ...\}}$ converges at $n = n_m$, one sets:
\begin{equation}
    \{\phi_{\alpha,m+1}^{0}, ...\} \coloneqq \{\phi_{\alpha,m}^{n_m}, ...\}
\end{equation}
and the new timestep $m+1$ starts.

\subsection{Algorithm}
\label{subsec:FEM_Algor}
This subsection demonstrates a more fine-grained algorithm that assembles multiple variables across different workpieces. 
For simplicity, we will still use the generalize-$\alpha$ temporal discretization 
and only show the process for the domain physics. 
In practice, the boundary conditions are just the domain physics one dimension lower and 
the numerical modifications are treated in the same way as the domain physics, so that 
the same process is essentially repeated another two times with different data.

\begin{figure}[tbhp]
    \centering
    \includegraphics[width=\textwidth]{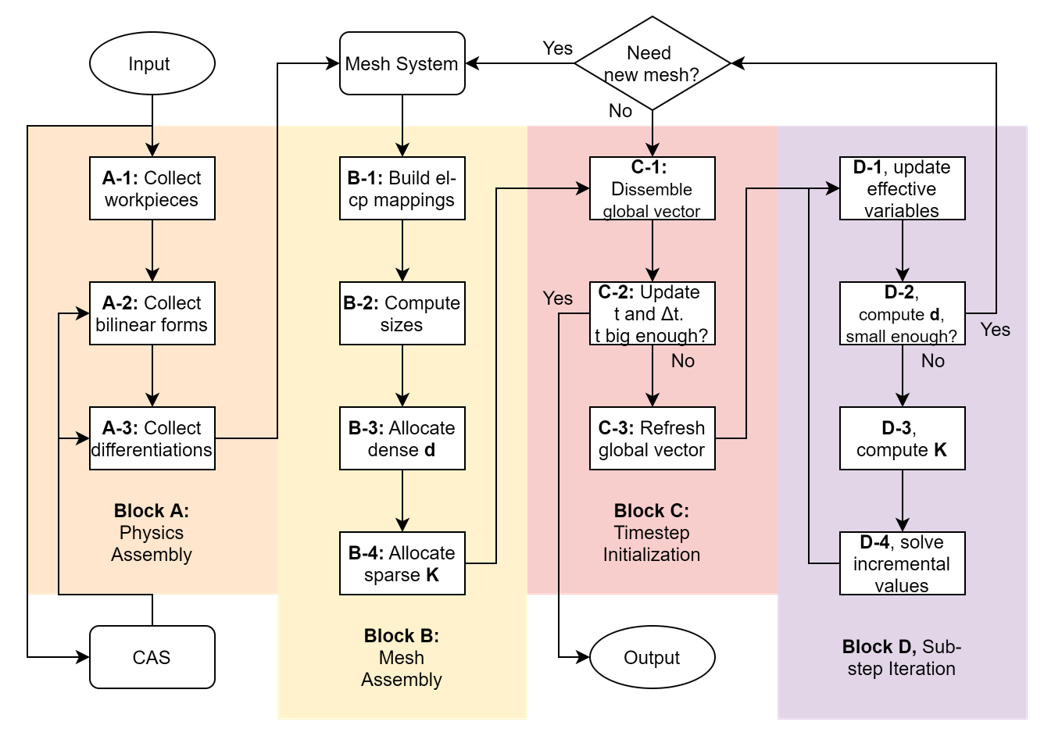}
    \caption{FEM kernel flowchart}
    \label{fig_FEM_kernel}
\end{figure}

The proposed FEM kernel consists of 4 Blocks, as shown in \cref{fig_FEM_kernel}, 
positioned from left to right and with abstraction levels from high to low as:
\begin{enumerate} 
    \item Physics Assembly;
    \item Mesh Assembly;
    \item Timestep Initialization; and
    \item Sub-step Iteration;
\end{enumerate}
which will be discussed one by one next.

\textbf{Block A} - \textit{Physics Assembly} is completely symbolic and only runs right after the physics is assigned to the workpieces 
(with or without discretization to the mesh). 

\textbf{A-1}, all workpieces are collected as $\{\Omega_{\beta^{wp}}\}$,  $\beta^{wp} =1,2,...,\hat{\beta}^{wp}$.

\textbf{A-2}, for each workpiece $\Omega_{\beta^{wp}}$, all the bilinear forms are collected as $\{B_{\beta^{b}}\}_{\beta^{wp}}$, 
where $\beta^{b}$ is not tracked since the bilinear forms will be re-organized by the rewriting rules.
Meanwhile, all the symbols of the basic variables are collected in $\{\phi^{\kappa}\}_{\beta^{wp}}$, 
with superscripted index $\kappa = 1,2,...,\hat{\kappa}_{\beta^{wp}}$ 
since the subscript is reserved for the control point index. 

\textbf{A-3}, for each bilinear form in $\{B_{\beta^{b}}\}_{\beta^{wp}}$, i.e.,
\begin{equation*}
(D_0\overline{\phi}^{\kappa_0}, \mathcal{L}^a(..., \partial_t^{\nu_{\lambda}}D_{\lambda}\phi^{\kappa_\lambda}, ...))_{\Omega_{\beta^{wp}}}      
\end{equation*}
where the $\lambda^{th}$ basic variable symbol index is specified by $\kappa_\lambda$, we
collect the unique pairs $\{(\kappa_0, \kappa_{\lambda})_{\beta^{sym\_pair}}\}_{\beta^{wp}}$ with index 
$\beta^{sym\_pair}=1,2,...,\hat{\beta}^{sym\_pair}_{\beta^{wp}}$, for example, by integer hashing, 
to assemble the sparse matrix $\mathbf{K}$ later.

The operator $\mathcal{L}^a$ is also differentiated with respect to each operand (word, not symbol), denoted by $\frac{\partial\mathcal{L}^a}{\partial\lambda}$, 
resulting in the bilinear forms: 
\begin{equation*}
    (D_0\overline{\phi}^{\kappa_0}, \frac{\partial\mathcal{L}^a}{\partial\lambda}(..., \partial_t^{\nu_{\lambda}}D_{\lambda}\phi^{\kappa_\lambda}, ...))
    _{\Omega_{\beta^{wp}}}      
\end{equation*}
which are collected in $\{B^{diff}_{\beta^{b'}}\}_{\beta^{wp}}$ for computing the matrix value $\mathbf{K}$ later.

\textbf{Block A} ends above. Before \textbf{Block B}, the extension to multiple variables and multiple workpieces can be briefly summarized as:
\begin{enumerate} 
    \item The residue $d$ is modified to the sum of residues from all workpieces:
    \begin{equation}
        d=\sum_{\beta^{wp} = 1}^{\hat{\beta}^{wp}} d_{\beta^{wp}}
    \end{equation}
    \item The global DOFs, denoted by $\alpha'$ and $\overline{\alpha}'$, are not only the workpiece control point indices $\alpha, \overline{\alpha}$, 
    but are determined by a mapping from the workpiece, basic variable and workpiece control point indices hierarchically as:
    \begin{equation}
        \alpha' \coloneqq \alpha'(\beta^{wp}, \kappa, \alpha) 
    \end{equation}
    \begin{equation}
        \overline{\alpha}' \coloneqq \overline{\alpha}'(\beta^{wp}, \overline{\kappa}, \overline{\alpha})
    \end{equation}
    \item The FEM is to minimize each discretized component in the global DOF, according to:
    \begin{equation}
        d_{\overline{\alpha}'(\beta^{wp}, \overline{\kappa}, \overline{\alpha})} \coloneqq 
        \frac{\partial{d_{\beta^{wp}}}}{\partial(\delta\overline{\phi}^{\overline{\kappa}}_{{\overline{\alpha}}})} 
    \end{equation}
\end{enumerate}

\textbf{Block B} - \textit{Mesh Assembly} runs after all the workpieces are discretized and whenever the mesh is changed.

\textbf{B-1}, for each workpiece $\Omega_{\beta^{wp}}$, 
\begin{enumerate}
    \item The element indices are $\beta^{el}=1,2, ..., \hat{\beta}^{el}_{\beta^{wp}}$;
    \item The workpiece control point indices are $\alpha=1,2,...,\hat{\alpha}_{\beta^{wp}}$;
    \item In each element $\beta^{el}$, the control points can be referred to by element to control point mapping 
    $\alpha=\alpha(\beta^{el}, \beta^{el\_cp})$ where $\beta^{el\_cp} = 1,2,...,\hat{\beta}^{el\_cp}_{\beta^{el}}$ 
    with $\hat{\beta}^{el\_cp}_{\beta^{el}}$ being the number of control points in this element; 
    \item The workpiece control point index pairs 
    $\{(\alpha_1, \alpha_2)_{\beta^{cp\_pair}}\}_{\beta^{wp}}$ are collected with $\beta^{cp\_pair} = 1,2,...,\hat{\beta}^{cp\_pair}_{\beta^{wp}}$ 
    by variating the last input in each element to control point mapping, 
    where each unique pair is only kept once, like the basic symbol pairs.
\end{enumerate}

\textbf{B-2}, For each workpiece $\Omega_{\beta^{wp}}$, the workpiece-wise last dense ID $n^{dense}_{\beta^{wp}}$ is: 
\begin{equation}
    n^{dense}_{\beta^{wp}}\coloneqq\sum_{\beta'=1}^{\beta^{wp}}\hat{\kappa}_{\beta'}\times{\hat{\alpha}_{\beta'}}
\end{equation}
and the workpiece-wise last sparse ID $n^{sp}_{\beta^{wp}}$ is: 
\begin{equation}
    n^{sp}_{\beta^{wp}}\coloneqq\sum_{\beta'=1}^{\beta^{wp}}\hat{\beta}^{sym\_pair}_{\beta'}\times{\hat{\beta}^{cp\_pair}_{\beta'}}
\end{equation}
Note that dense ID is for the residue vector $\mathbf{d}$, while the sparse matrix is for the stiffness matrix $\mathbf{K}$.

\textbf{B-3}, for the residue vector, the size is $\hat{\alpha}'= n^{dense}_{\hat{\beta}^{wp}}$ and the global DOF mapping is explicitly:
\begin{equation}
    \alpha'(\beta^{wp}, \kappa, \alpha) \coloneqq n^{dense}_{\beta^{wp}-1} + (\kappa-1)\hat{\alpha}_{\beta^{wp}}+\alpha
\end{equation}
\begin{equation}
        \alpha'(\beta^{wp}, \kappa, \beta^{el}, \beta^{el\_cp}) \coloneqq \alpha'(\beta^{wp}, \kappa, \alpha(\beta^{el}, \beta^{el\_cp}))
\end{equation}
where $n^{dense}_0 = 0$. 

The global effective arrays $\{\partial_t^{\nu}\widetilde{\phi}_{\alpha'}\}$, 
the increment variable arrays $\{\Delta\partial_t^{\nu}\phi_{\alpha'}\}$ and the global residue array $\{d_{\overline{\alpha}'}\}$,
are allocated with sizes $\hat{\nu}\times\hat{\alpha}'$, $\hat{\nu}\times\hat{\alpha}'$, and $\hat{\alpha}'$ respectively.
Note that here both $\partial_t^{\nu}\phi$ and $\partial_t^{\nu}\widetilde{\phi}$ are regarded as numerical variables instead of symbolic differentiations, i.e., 
$\{\partial_t\phi_{\alpha'}\} = \{u_{\alpha'}\}$, so as in the Block C and D.

\textbf{B-4}, for the stiffness matrix, we:
\begin{enumerate}
    \item Define the total sparse vector size $\hat{\beta}^{sp}=n^{sp}_{\hat{\beta}^{wp}}$;
    \item Allocate the row ID, the column ID and the value arrays $\{I_{\beta^{sp}}\}$, $\{J_{\beta^{sp}}\}$, $\{K_{\beta^{sp}}\}$ 
    with $\beta^{sp} = 1,2,...,\hat{\beta}^{sp}$ defining the COO-format sparse matrix;
    \item For each basic symbol pair and workpiece control point pair, i.e.,
    \begin{equation*}
        \{(\kappa_0, \kappa_{\lambda})_{\beta^{sym\_pair}}, \quad (\alpha_1, \alpha_2)_{\beta^{cp\_pair}}\}_{\beta^{wp}}
    \end{equation*}
    we define the sparse ID mapping
    \begin{equation}
        \beta^{sp}(\beta^{wp}, \beta^{sym\_pair}, \beta^{cp\_pair}) \coloneqq 
        n^{sp}_{(\beta^{wp}-1)} + (\beta^{sym\_pair} - 1)\hat{\beta}^{cp\_pair}_{\beta^{wp}} + \beta^{cp\_pair} \label{Eq:sparse_ID}
    \end{equation}
    and fill the corresponding component of the row ID array and the column ID array with:
    \begin{alignat}{4}
        I_{\beta^{sp}} =\alpha'(\beta^{wp}, \kappa_0, \alpha_1)={}& n^{dense}_{(\beta^{wp}-1)} &+& (\kappa_0 - 1)\hat{\alpha}_{\beta^{wp}} &+& \alpha_1 \\
        J_{\beta^{sp}} =\alpha'(\beta^{wp}, \kappa_\lambda, \alpha_2)={}& n^{dense}_{(\beta^{wp}-1)} &+& (\kappa_\lambda - 1)\hat{\alpha}_{\beta^{wp}} &+& \alpha_2 
    \end{alignat}
\end{enumerate}

\textbf{Block C} - \textit{Timestep Initialization} runs at the beginning of each timestep to refresh the incremental data.

\textbf{C-1}, each workpiece control point variable value $\partial_t^{\nu}\phi^{\kappa}_{\alpha}$ in $\Omega_{\beta^{wp}}$
is updated by:
\begin{equation}
    \partial_t^{\nu}\phi^{\kappa}_{\alpha}\mathrel{+}= \Delta\partial_t^{\nu}\phi_{\alpha'(\beta^{wp},\kappa,\alpha)}
\end{equation} 

\textbf{C-2}, the new time $t$ and timestep $\Delta{t}$ are calculated. $\Delta\partial_t^{\hat{\nu}}\phi_{\alpha'}$ is cleared (to zeros).

\textbf{C-3}, new $\{\Delta\partial_t^{\nu}\phi_{\alpha'}\}$ is initialized by updating the sequence of $\nu=(\hat{\nu} - 1), (\hat{\nu}-2),..,0$, i.e.: 
\begin{equation}
    \Delta\partial_t^{\nu}\phi_{\alpha'} = (b_{(\nu + 1)}\Delta{t})\partial_t^{(\nu + 1)}\phi^{\kappa}_{\alpha} + \Delta\partial_t^{(\nu + 1)}\phi_{\alpha'}
\end{equation}  

\textbf{Block D} - \textit{Sub-step Iteration} is executed at each sub-step.

\textbf{D-1}, the effective variable values $\{\partial_t^{\nu}\widetilde{\phi}_{\alpha'}\}$ are updated as:
\begin{equation}
    \partial_t^{\nu}\widetilde{\phi}_{\alpha'} = c_{(\nu + 1)}\Delta\partial_t^{\nu}\phi_{\alpha'} + \partial_t^{\nu}\phi^{\kappa}_{\alpha}
\end{equation}

\textbf{D-2}, the residue array $\{d_{\overline{\alpha}'}\}$ is cleared first. Then 
for $\overline{\beta}^{el\_cp} = 1,2,...,\hat{\beta}^{el\_cp}_{\beta^{el}}$ in each element $\beta^{el}$ of each bilinear form in $\{B_{\beta^{b}}\}_{\beta^{wp}}$
in each workpiece $\beta^{wp}$
\begin{equation*}
    (D_0\overline{\phi}^{\kappa_0}, \mathcal{L}^a(..., \partial_t^{\nu_{\lambda}}D_{\lambda}\phi^{\kappa_\lambda}, ...))_{\Omega_{\beta^{wp}}}      
\end{equation*}
the residue vector $\{d_{\overline{\alpha}'}\}$ is updated by the (atomic) increment:
\begin{align}
    &d_{\alpha'(\beta^{wp}, \kappa^0, \beta^{el}, \overline{\beta}^{el\_cp}))}\mathrel{+}= \\
    &\nonumber \sum_{\gamma}\{w_{\gamma}^{\text{itg}}(D_0\overline{N}_{\alpha(\beta^{el},\overline{\beta}^{el\_cp})})|_{\vec{x}_{\gamma}^{\text{itg}}}\\
    &\nonumber \mathcal{L}^a(...,\sum_{\beta^{el\_cp} = 1}^{\hat{\beta}^{el\_cp}_{\beta^{el}}}(D_\lambda N_{\alpha(\beta^{el},\beta^{el\_cp})})|_{\vec{x}_{\gamma}^{\text{itg}}}
    (\partial_t^{\nu_\lambda}\widetilde{\phi}_{\alpha'(\beta^{wp}, \kappa_{\lambda}, \beta^{el}, \beta^{el\_cp})}), ...)\}
\end{align} 
If the residue is small enough, break, else continue.

\textbf{D-3}, the sparse value array $\{K_{\beta^{sp}}\}$ is cleared first. 
Then for $\beta^{el\_cp}, \overline{\beta}^{el\_cp} = 1,2,...,\hat{\beta}^{el\_cp}_{\beta^{el}}$ in each element $\beta^{el}$ 
of each differentiated bilinear form in $\{B^{diff}_{\beta^{b'}}\}_{\beta^{wp}}$ in each workpiece $\beta^{wp}$,
\begin{equation*}
    (D_0\overline{\phi}^{\kappa_0}, \frac{\partial\mathcal{L}^a}{\partial\lambda^{diff}}(..., \partial_t^{\nu_{\lambda}}D_{\lambda}\phi^{\kappa_\lambda}, ...))_{\Omega_{\beta^{wp}}}      
\end{equation*}
the sparse value array component $K_{\beta^{sp}}$ is updated by the (atomic) increment:
\begin{align}
    & K_{\beta^{sp}}\mathrel{+}= \\
    &\nonumber \sum_{\gamma}\{w_{\gamma}^{\text{itg}}(D_0\overline{N}_{{\alpha(\beta^{el},\overline{\beta}^{el\_cp})}}
    D_{\lambda^{diff}}N_{\alpha(\beta^{el}, \beta^{el\_cp})})|_{\vec{x}_{\gamma}^{\text{itg}}} 
    \frac{c_{(\nu_{\lambda^{diff}}+1)}}{\prod_{\beta'= 1}^{\nu_{\lambda^{diff}}}(b_{\beta'}\Delta{t})} \times\\
    &\nonumber \frac{\partial\mathcal{L}^a}{\partial\lambda^{diff}}
    (...,\sum_{\beta'= 1}^{\hat{\beta}^{el\_cp}_{\beta^{el}}}(D_\lambda N_{\alpha(\beta^{el},\beta')})|_{\vec{x}_{\gamma}^{\text{itg}}}
    (\partial_t^{\nu_\lambda}\widetilde{\phi}_{\alpha'(\beta^{wp}, \kappa_{\lambda}, \beta^{el}, \beta')}), ...)\}
\end{align} 
where the sparse matrix entry $\beta^{sp}$ is computed from \cref{Eq:sparse_ID} with the basic symbol pair index $\beta^{sym\_pair}$ from pair $(\kappa_0, \kappa_{\lambda^{diff}})$ and
the local control point pair index $\beta^{cp\_pair}$ from pair $(\alpha(\beta^{el}, \overline{\beta}^{el\_cp}),\alpha(\beta^{el}, \beta^{el\_cp}))$.

\textbf{D-4}, the sub-step increment $\{\Delta_{sub}\phi_{\alpha'}\}$ is calculated by solving $\mathbf{K}/\mathbf{d}$ through a (direct or iterative) linear solver.
The incremental variable array $\{\Delta\partial_t^{\nu}\phi_{\alpha'}\}$ is updated by:
\begin{equation}
    \Delta\partial_t^{\nu}\phi_{\alpha'} \mathrel{+}= \Delta_{sub}\phi_{\alpha'}\prod_{\beta'= 1}^{\nu}(b_{\beta'}\Delta{t}),\quad\nu = 0,1,...,\hat{\nu}
\end{equation}
and the program returns to the Block start \textbf{D-1}.

With the four (4) blocks discussed above, we can have a brief review on the kernel workflow in \cref{fig_FEM_kernel}. \textbf{Block A}
parses and rewrites (automatically derives) the input expressions of the physics into intermediate representations. 
\textbf{Block B} links the intermediate representations
to the mesh and allocates all the data arrays. \textbf{Block C} updates the data arrays by one timestep by repeating 
\textbf{Block D} until the residue is reduced below the tolerance. Or more function-wise, firstly \textbf{Block A} processes symbols,
then \textbf{Block B} assembles arrays, finally \textbf{Block C} and \textbf{Block D} actually run the FEM simulation. 
Now, we will discuss in detail in \cref{sec:RW_system} how exactly symbols are processed.

\section{The rewriting system}
\label{sec:RW_system}

Besides of parsing the input mathematical expressions at the beginning and generating the corresponding codes at the end,
the main symbolic functions of the rewriting system are simplification (in multiple ways) and symbolic differentiation.
This section demonstrates a rewriting system designed primarily from a continuum mechanics perspective following the notations
used in \cite{RW_system_1999}.

\subsection{Theory} 
The fundamental elements of the rewriting systems are symbols with (or without) indices, denoted by words.
A word can be represented in the meta form $w_{a_1...a_n,b_1...b_m}$, where $w$ is the physical tensor symbol, 
$a_1...a_n$ are the component indices with $n$ being the tensor order and $b_1...b_m$ are the derivatives.
As the data structure implies, covariance and contravariance indices are not distinguished 
in MetaFEM for simplicity.

Next, the signature $\Sigma$ is defined by the union of:
\begin{enumerate}
    \item The set of all words, denoted by $\mathcal{W}$ (W for word);
    \item The set of all real (floating-point) numbers, denoted by $\mathcal{N}$ (N for number); and
    \item The set of all function names, denoted by $\mathcal{F}$ (F for function), which contains $\mathcal{L}^a$ in \cref{sec:FEM_kernel}.
\end{enumerate}

A ground term $T$ is defined to be either an element in $\mathcal{W}\cup\mathcal{N}$, or a tree with a function name in $\mathcal{F}$
as the tree root and other ground terms being the leaves,
which may result in a nested tree structure. The collection of ground terms is denoted as $\mathbb{T}(\Sigma, \emptyset)$.

\begin{figure}[tbhp]
    \centering
    \includegraphics[width=0.4\textwidth]{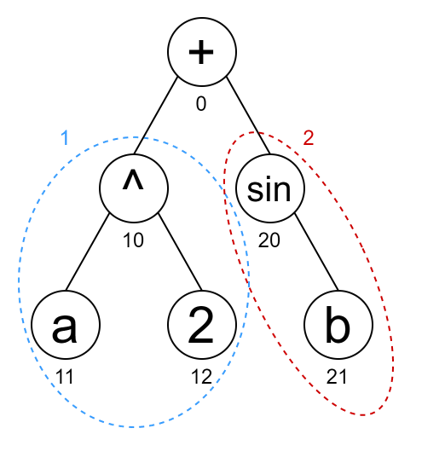}
    \caption{The example ground term $T=a^2+\sin(b)$ with positions.}
    \label{fig_tree}
\end{figure}

An example ground term $T=a^2+\sin(b)$ is shown in \cref{fig_tree}.
Each node of $T$ can be referred to by a sequence of indices, denoted by a position. 
The set of all the possible positions is denoted by $Pos(T)=\{0,1,2,10,11,12,20,21\}$, $T|_1=a^2$, $T|_{20}=\sin$, $T|_{21}=b$, etc., where index $0$ 
represents the function name while a positive integer index $n$ represents the $n^{th}$ subnode. 
Note, the separator in the sequence of indices is omitted when the maximum subnode number is under 10,
following the classical literature, 
so that $10,11,...$ in the above $Pos(T)$ are actually "$1$-$0,1$-$1,...$".

In practice, a ground term may represent an algebraic operation $\mathcal{L}^a(...)$, 
a bilinear form $(\cdot,\cdot)$, the overall weak form (but without the geometry part), or more generally,
a term with some concrete physical meaning a program can evaluate. In contrast, a general term is 
a term optionally with some placeholders to provide structure information, but cannot be evaluated until all the placeholders are filled, 
in which case the general term becomes a ground term.

The placeholders are named as syntactical variables $\mathsf{X}$, 
which can be further divided into placeholders for ground terms $\mathsf{W}$ and placeholders for operators $\mathsf{F}$, 
$\mathsf{X}=\mathsf{W}\cup\mathsf{F}$ (new fonts are used for new definitions).
With the syntactical variables, the general terms 
$\mathbb{T}(\Sigma, \mathsf{X})$ are defined to be the set where 
each general term $T\in\mathbb{T}(\Sigma, \mathsf{X})$ is either an element in 
$\mathcal{W}\cup\mathcal{N}\cup\mathsf{W}$ or a tree with an element in $\mathcal{F}\cup\mathsf{F}$ being the tree root and other general terms being the leaves. 

A ground term is a general term with no variables, therefore, the general term $T$ can also be referred to by the positions in the same way.
With the positions $Pos(T)$, the set of all occurring variables of term is defined by:
\begin{equation}
    Var(T)\coloneqq\{v|\quad\exists p\in{Pos(T)},\quad v=T|_p,\quad v\in\mathsf{X}\}
\end{equation}

On the other hand, a general term can be filled to generate a new ground term. 
The generalized filling process is denoted by substitution $\sigma=\sigma_{(\mathsf{X}^0,\mathcal{X}^0)}$
which is determined by a variable list 
$\mathsf{X}^0=\{W_1,...,W_m, F_1, ..., F_n\}$ and a result list 
$\mathcal{X}^0=\{T_1,...,T_m, G_1, ..., G_n\}$, 
where $W_1,...,W_m\in\mathsf{W}$, $F_1,...,F_n\in\mathsf{F}$, $T_1,...,T_m\in\mathbb{T}(\Sigma, \emptyset)$ and
$G_1,...,G_n\in\mathcal{F}$ with $m,n$ being non-negative integers. 

The substitution can operate on a term by:
\begin{alignat}{4}
    \nonumber &\sigma(W_p)&{}={}&T_p&p={}& 1,...,m\\
    \nonumber &\sigma(F_q)&{}={}&G_q&q={}&1,...,n\\
    \nonumber &\sigma(V)&{}={}&V&V\in{}&(\Sigma\cup\mathsf{X})-\mathsf{X}^0\\
    \sigma(F(V_1,V_2&,...))&{}={}&\sigma(F)&(\sigma(V_1),\sigma(V_2),...))&
\end{alignat}
In other words: (1) when $\sigma$ acts on a term tree, it keeps the tree structure and acts on each node; and (2) when $\sigma$ 
acts on a variable, if the variable is found in $\mathsf{X}^0$, $\sigma$ returns the element at the same index in $\mathcal{X}^0$, 
otherwise it does nothing.

Next, the match operator $M$ can be defined to determine whether a general term $T^M$ represents a ground term $T$,
that is, $M(T^M, T)=(true, \mathcal{X}^M)$ if there exists a mapping list $\mathcal{X}^M$ so that:
\begin{equation}
    \sigma_{(Var(T^M),\mathcal{X}^M)}(T^M)=T
\end{equation}
otherwise $M(T^M, T)=(false, \emptyset)$. 
Note $\mathcal{X}^M$ can also be empty like in $M(T^M, T)=(true,\emptyset)$ when $T^M$ is a specific ground term.

\begin{figure}[tbhp]
    \centering
    \includegraphics[width=0.4\textwidth]{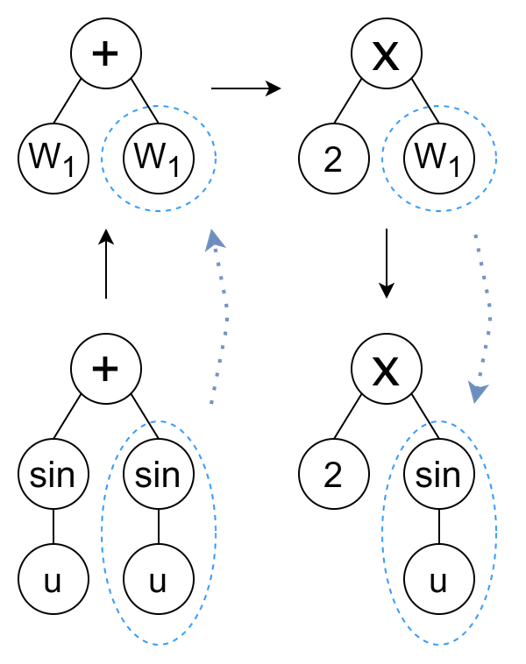}
    \caption{The example rule rewriting $sin(u) + sin(u)\rightarrow2 * sin(u)$.}
    \label{fig:CAS_rule}
\end{figure}

A classical rewrite rule $R$ is determined by a suitable pair of general terms $T^l, T^r$, 
where the superscripts are named for the Left Hand Side (LHS) and Right Hand side (RHS) in an equation, 
although the considered rewrite rules are one-sided derivations (left to right only) instead of real equations,
as shown in \cref{fig:CAS_rule},

where:
\begin{alignat}{2}
    \nonumber R\coloneqq(T^l, T^r)\in\mathbb{T}(\Sigma,\mathsf{X})\times\mathbb{T}(\Sigma,\mathsf{X})\\
    Var(T^r)\subseteq Var(T^l)
    \label{Eq:varset}
\end{alignat}
and $R$ can operate on a ground term node $T$ (but not its subnodes) by:
\begin{align}
    \nonumber (Bool, \mathcal{X}^l)\coloneqq{}&M(T^l, T)\\
    R(T)={}&
    \begin{cases}
        \sigma_{(Var(T^l),\mathcal{X}^l)}(T^r),&Bool=true\\
        T,&Bool=false
    \end{cases}
\end{align}
which is, if $T$ matches $T^l$, $R$ rewrites $T$ by substitution of $T^r$ with the mapping list from $M(T^l, T)$, 
otherwise it does nothing.  

Further, a rewrite rule $R$ can recursively apply to a whole ground term $T$, which is denoted by $\widehat{R}(T)$, 
e.g., a Leftmost-Innermost (LI) formulation is:
\begin{equation}
    \widehat{R}(T)=
    \begin{cases}
        R(F(\widehat{R}(V_1),\widehat{R}(V_2),...)),&T=F(V_1, V_2,...)\\
        R(T),&\text{else}
    \end{cases}
\end{equation} 
A practical rewriting is to repetitively apply a set of rewriting rules to the target term until no change occurs.

In practice, one may (plainly) extend the classical rewrite rule $(T^l, T^r)$ in two ways:
\begin{enumerate}
    \item The match operator $M$ is purely about structure. An operator can be defined to check the content (per rewrite rule), 
    denoted by the checker $Ck\coloneqq\{c_1, c_2,...,c_{m+n}\}$, where each $c_{\beta}, \quad\beta=1,2,...,m+n$ is simply a predefined set.
    The checker checks whether each entry in the (matched) LHS result list $\mathcal{X}^l=\{T^l_1,...,T^l_m, G^l_1,...,G^l_n\}$ belongs to the set 
    at the same index of itself as follows:
    \begin{equation}
        Ck(\mathcal{X}^0)=
        \begin{cases}
            true,&T^l_1\in c_1, ..., T^l_m\in c_m, G^l_1\in c_{m+1}, ..., G^l_n\in c_{m+n}\\
            false,&\text{else},
        \end{cases}
    \end{equation}
    \item The restriction in \cref{Eq:varset} keeps $R(T)$ a ground term, but is too strict. What is really needed is that $Var(T^r)$ is computable from $Var(T^l)$.
    Therefore, for a general $Var(T^r)=\{W^r_1, ..., W^r_{m^r}, F^r_1, ..., F^r_{n^r}\}$, 
    \cref{Eq:varset} can be replaced by a operator (per rewrite rule), denoted by the transcriber $Tb\coloneqq\{H_1,...,H_{m^r+n^r}\}$ with each entry a (programmable) 
    function which takes in the LHS result list $\mathcal{X}^l$ and outputs the ground term at the same index to form the RHS result list $\mathcal{X}^r$, that is:
    \begin{equation}
        \mathcal{X}^r(\mathcal{X}^l)=Tb(\mathcal{X}^l)=\{H_1(\mathcal{X}^l),...,H_{m^r+n^r}(\mathcal{X}^l)\}
    \end{equation}
\end{enumerate}

Finally, the extended rewrite rule for the FEM CAS can be redefined by $R=(T^l, T^r, Ck, Tb)$ so that:
\begin{align}
    \nonumber (Bool,\mathcal{X}^l)\coloneqq{}&M(T^l, T)\\
    R(T)={}&
    \begin{cases}
        \sigma_{(Var(T^r),Tb(\mathcal{X}^l))}(T^r),&Bool\land Ck(\mathcal{X}^l)=true\\
        T,&\text{else}
    \end{cases}
\end{align}
which is ready for implementation with unchanged $\widehat{R}(T)$. A brief review of the proposed rewriting system is listed 
in \cref{tab:RW_system}.

\begin{table}[tbhp]
    {\footnotesize
      \caption{Overview of the rewriting system.} \label{tab:RW_system}
    \begin{center}
      \begin{tabular}{|c|c|c|c|c|c|} \hline
        \multicolumn{3}{|c}{Signature} & \multicolumn{2}{|c|}{Terms} & Rules \\ \hline 
        Numbers & Words & Functions & Ground terms & General terms & \\ \hline
        $0.65, 1.0$ & $p,u_{i,t},$ & $+,\times,(\cdot,\cdot),$ & $0.1 + p,$ & $a^2+sin(b),$ & $W_1+0\rightarrow W_1,$ \\ 
        &$\sigma_{ij,j}$& $my func$ & $a^2+sin(b)$ & $W_1, W_1 + W_1$ & $ W_1 + W_1\rightarrow 2\times W_1$ \\ \hline
      \end{tabular}
    \end{center}
    }
\end{table}

\subsection{Implementation}
The code for a rewriting system can be very simple and straightforward, 
since applying a rewriting rule $R(\cdot)$ is essentially a pattern matching problem from the computing perspective.
An example data structure, instead of the real code to make it simpler in few lines, is shown below.

There is no difference in data structure between placeholders for words or for functions, and only a general variable is needed.
On the other hand, although the checker and the transcriber $Ck, Tb$ are defined as rule-wise collections, each $c$ and each $H$ 
are more conveniently attached to the corresponding variables in practice. Therefore, a syntactic variable is defined by:
\begin{jllisting}
    mutable struct Syntactic_Variable
        id::Symbol
        c::Function
        H::Function
        H_args::Vector{Syntactic_Variable}
    end
\end{jllisting}
where $H\_args$ contains the arguments of $H$.

A term is either a number, a word (definition \cite{MetaFEM} omitted here for simplicity), a variable or a tree, and can be formulated as:
\begin{jllisting}
    struct Term_Tree
        head::Union{Function_Name, Syntactic_Variable}
        args::Vector{Union{Term_Tree, Number, Word, Syntactic_Variable}}
    end
    Term = Union{Term_Tree, Number, Word, Syntactic_Variable}
\end{jllisting}

Since $c$ and $H$ are variable-wise, a rewrite rule is simply:
\begin{jllisting}
    struct Rewrite_Rule
        LHS::Term
        RHS::Term
        LHS_vars::Vector{Syntactic_Variable}
        RHS_vars::Vector{Syntactic_Variable}
    end
\end{jllisting}
where $LHS\_vars, RHS\_vars$ are explicitly collected for illustration convenience since they contain both the checkers and 
transcribers. 

Locally applying a rewrite rule to a ground term node (but not its subnodes), i.e., $R(T)$, can be simply:
\begin{jllisting}
    function apply(rule::Rewrite_Rule, src_term::Term) 
        match_is_successful, LHS_result_list = match(rule.LHS, src_term)
        if match_is_successful
            checks_are_passed = check(LHS_result_list, rule.LHS_vars) 
            if checks_are_passed
                RHS_result_list = transcribe(LHS_result_list, rule.RHS_vars)
                return substitute(rule.RHS, RHS_result_list)
            end 
        end
        return src
    end 
\end{jllisting}
which can be described in 4 steps:
\begin{enumerate}
    \item Matching the rewrite rule LHS to the ground term node-by-node through a Deterministic Finite Automaton (DFA); 
    \item If the matching (of structures) is successful, check the LHS result list with checkers in LHS variables;
    \item If all the checks pass, transcribe the RHS result list by the transcriber in RHS variables and the LHS result list; and
    \item Generate the new ground term by substituting the rewrite rule RHS with the RHS result list.
\end{enumerate} 

With the local apply function, the global apply ($\widehat{R}(T)$) can be coded with multiple dispatch as:
\begin{jllisting}
apply_g(rule::Rewrite_Rule, src::Union{Word, Number}) = apply(rule, src)
apply_g(rule::Rewrite_Rule, src::Term_Tree) = apply(rule, 
    Term_Tree(src.head, map(x -> apply_g(rule, x), src.args))) 
\end{jllisting}

Finally, an intuitive grammar is designed for inputting rules in practice, and 
we show a real code that uses the current approach in closing this section. Specifically:
\begin{enumerate}
    \item (Part of) the code for expression simplification is:
    \begin{jllisting}
    @Define_Semantic_Constraint N₁ ∊ (N₁ isa Number)
    @Define_Semantic_Constraint N₂ ∊ (N₂ isa Number)

    @Define_Aux_Semantics N₁₊₂(N₁,N₂) = N₁ + N₂
    @Define_Aux_Semantics N₁₊₊(N₁) = N₁ + 1
    @Define_Aux_Semantics N₁₂(N₁,N₂) = N₁ * N₂
    @Define_Aux_Semantics N₁ₚ₂(N₁,N₂) = N₁ ^ N₂

    @Define_Rewrite_Rule Unary_Sub ≔ -(a) => (-1) * a
    @Define_Rewrite_Rule Binary_Sub ≔ (b - a => b + (-1) * a)
    @Define_Rewrite_Rule Binary_Div ≔ (a / b => a * b ^ (-1))

    @Define_Rewrite_Rule Unary_Add ≔ (+(a) => a)
    @Define_Rewrite_Rule Add_Splat ≔ ((a...) + (+(b...)) + (c...) => a + b + c)
    @Define_Rewrite_Rule Add_Numbers ≔ ((a...) + N₁ + (b...) + N₂ + (c...) => 
                                        N₁₊₂ + a + b + c)
    @Define_Rewrite_Rule Add_Sort ≔ (a + (b...) + N₁ + (c...) => N₁ + a + b + c)
    @Define_Rewrite_Rule Add_0 ≔ (0 + (a...) => +(a))

    @Define_Rewrite_Rule Unary_Mul ≔ (*(a) => a)
    @Define_Rewrite_Rule Mul_Splat ≔ ((a...) * (*(b...)) * (c...) => a * b * c)
    @Define_Rewrite_Rule Mul_Numbers ≔ ((a...) * N₁ * (b...) * N₂ * (c...) => 
                                        N₁₂ * a * b * c)
    @Define_Rewrite_Rule Mul_Sort ≔ (a * (b...) * N₁ * (c...) => N₁ * a * b * c)
    @Define_Rewrite_Rule Mul_0 ≔ (0 * (a...) => 0)
    @Define_Rewrite_Rule Mul_1 ≔ (1 * (a...) => *(a))

    @Define_Rewrite_Rule Pow_a0 ≔ (a ^ 0 => 1)
    @Define_Rewrite_Rule Pow_a1 ≔ (a ^ 1 => a)
    @Define_Rewrite_Rule Pow_0a ≔ (0 ^ a => 0)
    @Define_Rewrite_Rule Pow_1a ≔ (1 ^ a => 1)
    @Define_Rewrite_Rule Pow_Numbers ≔ (N₁ ^ N₂ => N₁ₚ₂)
    @Define_Rewrite_Rule Pow_Splat ≔ ((a ^ b) ^ c => a ^ (b * c))
    \end{jllisting}
    \item The code for symbolic differentiation is:
    \begin{jllisting}
    @Define_Semantic_Constraint S₁ ∊ (S₁ isa SymbolicWord)
    @Define_Semantic_Constraint S₂ ∊ (S₂ isa SymbolicWord)
    
    @Define_Semantic_Constraint f_undef ∊ (~(f_undef in [:+; :-; :*;
                                             :/; :^; :log; :inv]))
    
    @Define_Rewrite_Rule Basic_Diff1 ≔ ∂(S₁, S₁) => 1
    @Define_Rewrite_Rule Basic_Diff2 ≔ ∂(S₂, S₁) => 0
    @Define_Rewrite_Rule Basic_Diff3 ≔ ∂(N₁, S₁) => 0
    
    @Define_Rewrite_Rule Basic_Diff4 ≔ ∂({f_undef}(a...), S₁) => 0
    
    @Define_Rewrite_Rule Add_Diff ≔ ∂(a + (b...), S₁) => 
        ∂(a, S₁) + ∂(+(b...), S₁)
    @Define_Rewrite_Rule Mul_Diff ≔ ∂(a * (b...), S₁) => 
        ∂(a, S₁) * (b...) + a * ∂(*(b...), S₁)
    @Define_Rewrite_Rule Pow_Diff ≔ ∂(a ^ b, S₁) => 
        ∂(a, S₁) * a ^ (b - 1) * b  + ∂(b, S₁) * log(a) * a ^ b
    \end{jllisting}
\end{enumerate} 
where:
\begin{enumerate}
    \item Each line with \textit{Define Semantic Constraint} attaches the check function $c$ to the variable; 
    \item Each line with \textit{Define Aux Semantics} attaches the transcribe function $H$ and the function 
    argument variable $H\_args$ to the variable; and
    \item Each line with \textit{Define Rewrite Rule} defines a rule with the rule name before $\coloneqq$ 
    while the LHS and the RHS are connected with $\Rightarrow$ after $\coloneqq$. The splatting operator (...) denotes a variable quantity (0 or arbitrarily more) of function inputs, 
    similar to the Kleene star in a regular expression. 
\end{enumerate}

\section{Numerical results}
\label{sec:results}
The code is open-sourced on Github \cite{MetaFEM} with the modules discussed above and an original generic mesh system 
implemented with 2D/3D Lagrange cube/simplex elements of arbitrary order and serendipity cube elements 
of order 2 and 3 (but no mesh generation). 

All the six example cases shown later are also provided with most relevant files, including:
\begin{enumerate}
    \item Mesh and script for MetaFEM (this work) simulations;
    \item Commercial software setup and result files for comparison cases; and
    \item VTK files and Paraview states for visualization
\end{enumerate}
to be repeated in exact details. To be concise, we will only focus on the input code about physics in the ensuing discussions.

\subsection{Thermal conduction in a solid}
We begin with thermal conduction in a solid since it is among the simplest and most stable types of practical physics for FEM. 
The mathematical formulation is:
\begin{equation}
    \nonumber \text{variable}\quad T,\qquad\text{parameters}\quad C,k,h,h_{penalty},s,e_m,T_{env},T_{fix}
\end{equation}
\begin{alignat}{4}
    \nonumber -C(T,T_{,t})-k(T_{,i},T_{,i})+(T,s)={}&0,\qquad in\quad\Omega\\
    \nonumber h(T,T_{env}-T)+e_m\sigma^b(T, T_{env}^4 - T^4)={}&0,\qquad on\quad\partial(\Omega)_{convection\_radiation}\\
    h_{penalty}(T,T_{fix}-T)+k(T,n_iT_{,i})={}&0,\qquad on\quad\partial(\Omega)_{fix}
\end{alignat}
where $T$ is the temperature, $C$ is the volumetric heat capacity, $k$ is the thermal conductivity, $h$ is the convective
coefficient, $h_{penalty}$ is a large number to enforce the fixed temperature (Dirichlet) boundary condition, $s$ is the heat source,
$e_m$ is the emissivity and $\sigma^b=5.670\times10^{-8}W/m^{2}K^{4}$ is the Stefan–Boltzmann constant.

The code for the physics is:
\begin{jllisting}
@Def begin
    heat_dissipation = - C * Bilinear(T, T{;t}) - k * Bilinear(T{;i}, T{;i}) + 
                         Bilinear(T, s)
    conv_rad_boundary = h * Bilinear(T, Tₑₙᵥ - T) + em * σᵇ * Bilinear(T, Tₑₙᵥ^4 - T^4)
    fix_boundary = h_penalty * Bilinear(T, Tw - T) + k * Bilinear(T, n{i} * T{;i}) 
end
\end{jllisting}
Apparently, the grammar is almost exactly the same with the corresponding math, 
where the function Bilinear denotes a bilinear form $(\cdot,\cdot)$,
$\{\}$ in the code denotes subscripts in the math with the comma "," as a separator of words and 
the semicolon ";" marks the start of derivative indices. 

\subsubsection{Various boundary conditions on a 2D stripe}
A $2\,cm\times1\,cm$ rectangular strip is simulated for equilibrium with the bottom side insulated, 
the left, right sides fixed at $(900+273.15)\,K$ and the top side 
under convection and radiation with $T_{env}=(50+273.15)\,K,\quad h=50\,W/m^2K$ 
and emissivity $e_m=0.7$. 

\begin{figure}[tbhp]
    \centering
    \includegraphics[width=0.9\textwidth]{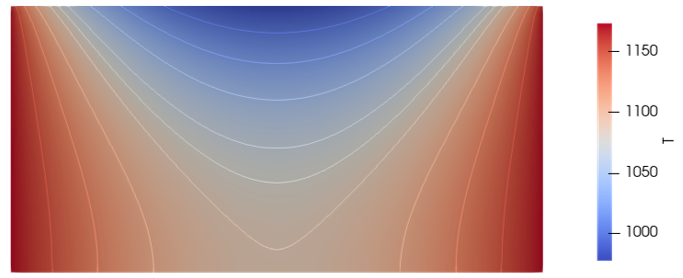}
    \caption{$T$ contour of the stripe.}
    \label{fig:2D_thermal_contour}
\end{figure}

\begin{figure}[tbhp]
    \centering
    \includegraphics[width=0.9\textwidth]{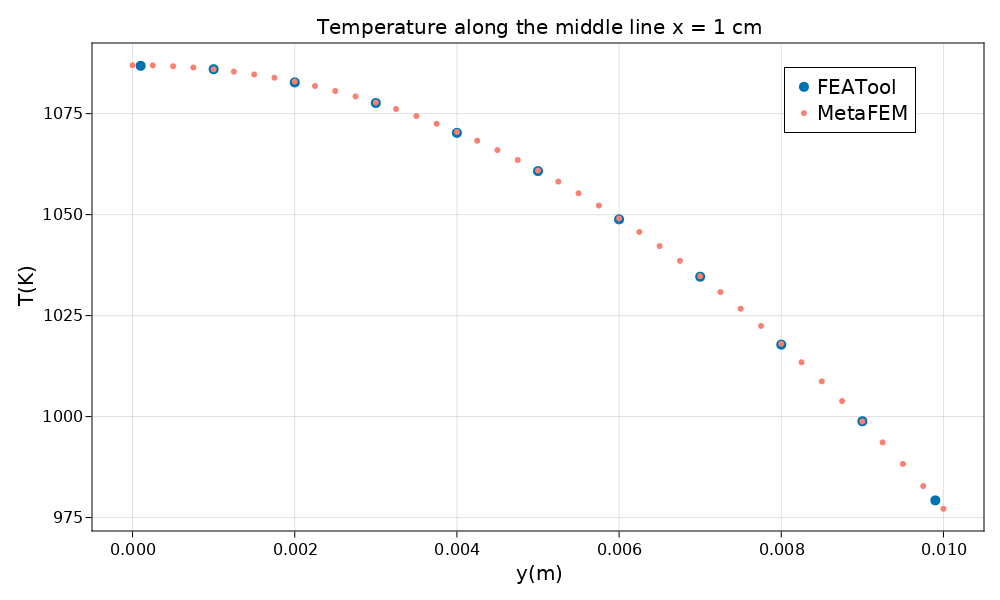}
    \caption{$T$ distribution along the vertical mid-line of the stripe. FEATool data are at sample points while MetaFEM (this work) data are at mesh nodes.}
    \label{fig:2D_thermal_lines}
\end{figure}

The case is a tutorial example in the matlab package FEATool \cite{FEATool_thermal}. 
Solved in MetaFEM by $40\times20$ quadratic serendipity rectangular elements (dx = $0.5\,mm$),
the temperature contour is plotted in \cref{fig:2D_thermal_contour} and 
the temperature distribution along the vertical mid-line (x = $1\,cm$) is 
compared to the FEATool result in \cref{fig:2D_thermal_lines}. The good match in temperature comparison indicates that
the domain physics and boundary conditions are correctly derived and assembled, including the nonlinear radiation term.

\subsubsection{Irregular 3D geometry with heat source}
To test the unstructured 3D mesh system, a pikachu is simulated where the body is considered as a uniform heat source $s=1.6\times10^3\,W/m^3$, $k=0.6\,W/(m\cdot K)$
from his own metabolism while the whole surface is under convection with $T_{env}=(20+273.15)\,K$, $h=25\,W/m^2K$. The numbers of property parameters are chosen 
by reference to the property parameters of a human in air.

\begin{figure}[tbhp]
    \centering
    \includegraphics[width=0.5\textwidth]{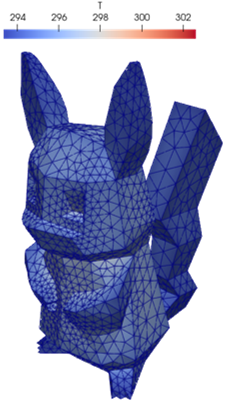}
    \caption{The mesh of pikachu, colored with temperature.}
    \label{fig:3D_thermal_mesh}
\end{figure}

\begin{figure}[tbhp]
    \centering
    \includegraphics[width=\textwidth]{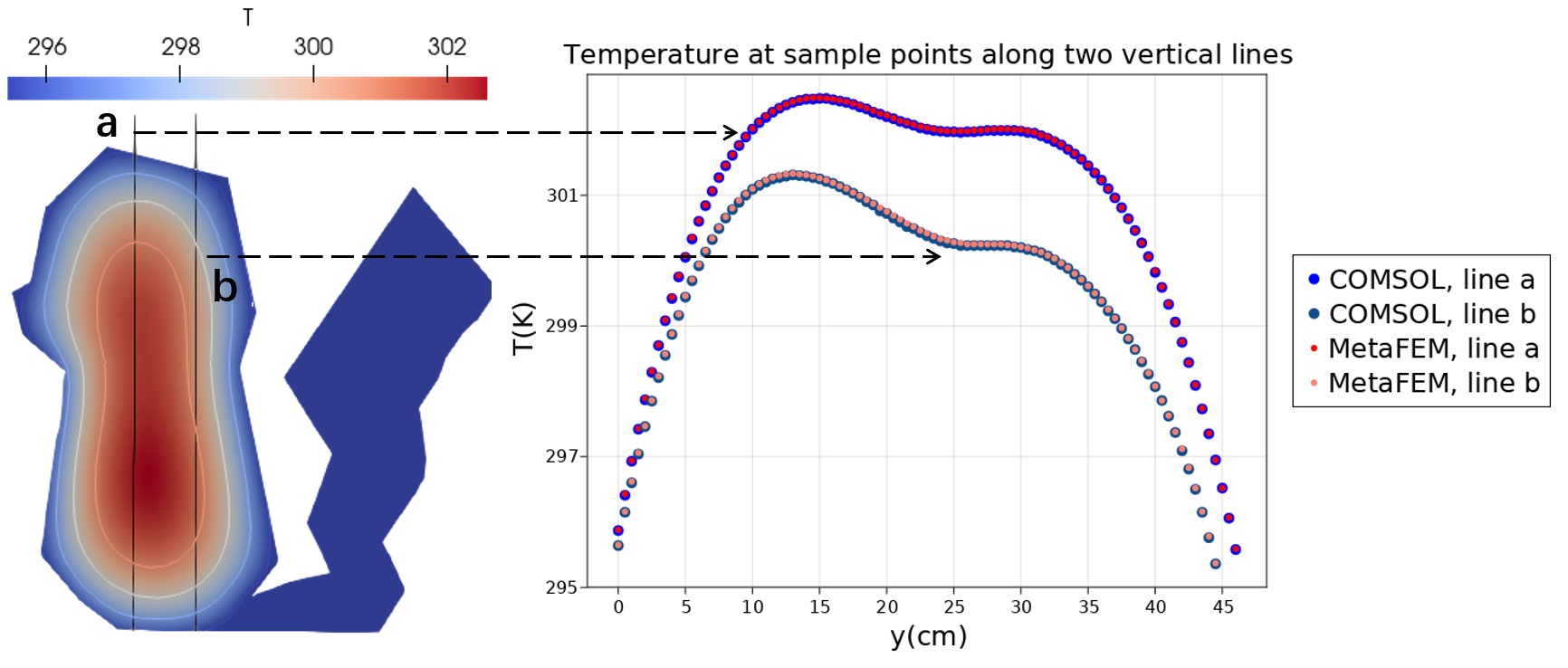}
    \caption{Pikachu temperature distribution. Left: Temperature contour on a slice; 
    Right: $T$ comparison at sample points along two vertical lines a and b in the pikachu.}
    \label{fig:3D_thermal_lines}
\end{figure}

The original source is a low-poly pikachu CAD file \cite{Pika_CAD} which is imported in COMSOL \cite{COMSOL_web} for both mesh generation and simulation. Then,
MetaFEM reads the same mesh and simulates with the same parameters. The mesh has 15,334 quadratic simplex elements and 23,703 nodes, 
which is shown in \cref{fig:3D_thermal_mesh} while the temperature distribution along sample 
lines a and b are compared across the two results in \cref{fig:3D_thermal_lines}. The good match indicates that the mesh system is able to handle complex, 
practical geometries equally well as common commercial solvers.

\subsection{Linear elasticity}
The mathematical formulation for the linear elastostatics in FEM is as follows:
\begin{equation}
    \nonumber \text{variable}\quad d_i,\qquad\text{parameters}\quad E,\nu,d^w_i,\sigma^l_{ij},\tau
\end{equation}
\begin{alignat}{4}
    \nonumber \lambda={}&\frac{E\nu}{(1+\nu)(1-2\nu)},\qquad&\mu={}&\frac{E}{2(1+\nu)}\\
    \nonumber \epsilon_{ij}={}&\frac{d_{i,j}+d_{j,i}}{2},\qquad&\sigma_{ij}={}&\lambda\delta_{ij}\epsilon_{mm}+2\mu\epsilon_{ij}
\end{alignat}
\begin{alignat}{4}
    \nonumber -(d_{i,j},\sigma_{ij})={}&0,\qquad in\quad\Omega\\
    \nonumber (d_i,\sigma^l_{ij}n_j)={}&0,\qquad on\quad(\partial\Omega)_{load}\\
    (d_i,\tau(d^w_i-d_i))={}&0,\qquad on\quad(\partial\Omega)_{fix}
\end{alignat}
where $d_i,E,\nu,\lambda,\mu,\epsilon_{ij},\sigma_{ij}$ are the displacement, Young’s modulus, Poisson ratio, Lame's first parameter, shear modulus, strain 
and stress, respectively.

The code for the physics is:
\begin{jllisting}
λ = E * ν / ((1 + ν) * (1 - 2 * ν))
μ = E / (2 * (1 + ν))

@Def ε{i,j} = (d{i;j} + d{j;i}) / 2.
@Def σ{i,j} = λ * δ{i,j} * ε{m,m} + 2. * μ * ε{i,j}

@Def begin
    ES_domain = - Bilinear(ε{i,j}, σ{i,j})
    ES_total_d_fixed_bdy = τ * Bilinear(d{i}, (dʷ{i} - d{i}))
    ES_d1_fixed_bdy = τ * Bilinear(d{1}, (dʷ{1} - d{1}))
    ES_loaded_bdy = Bilinear(d{i}, σˡ{i,j} * n{j})
end
\end{jllisting}

\subsubsection{3D cantilever beam bending}
\begin{figure}[tbhp]
    \centering
    \includegraphics[width=\textwidth]{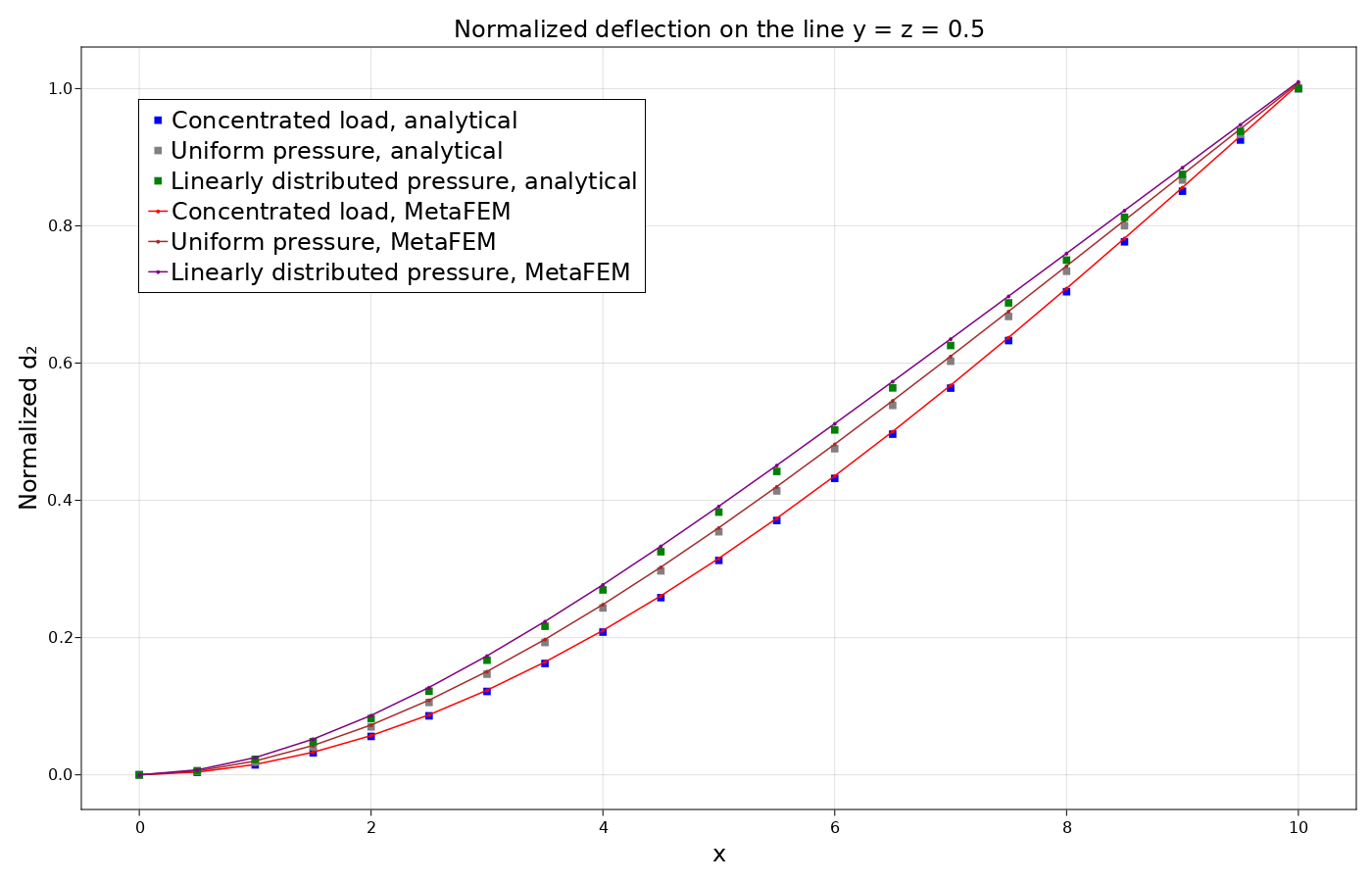}
    \caption{Cantilever bending. The analytical and numerical deflection distributions along the center line on the beam mid-plane after normalizing by 
    the corresponding maximum analytical deflection (to fit different deflection scales in one graph).}
    \label{fig:2D_Bending}
\end{figure}

A beam of normalized/dimensionless length $L=10$, width and depth $h=1$ with its left side (displacement) fixed is simulated for its deflection distribution along the mid-plane $y=z=0.5$
with $10\times4\times4$ quadratic serendipity brick elements under 3 different loading conditions:
\begin{enumerate}
    \item Concentrated load $P$ on the right side, where the analytical deflection distribution from beam theory \cite{beam_plate_1976} is:
    \begin{equation}
        d^{a1}_2(x)=-\frac{Px^2}{6EI}(3L-x),\qquad I\coloneqq\frac{1}{12h^3}
    \end{equation}
    \item Uniform pressure $p$ on the top, where the analytical deflection distribution is:
    \begin{equation}
        d^{a2}_2(x)=-\frac{px^2}{24EI}(6L^2-4Lx+x^2)
    \end{equation}
    \item Linear pressure $p'(x)=p_0(1-\frac{x}{L})$ on the top, where the analytical deflection distribution is:
    \begin{equation}
        d^{a3}_2(x)=-\frac{p_0x^2}{120LEI}(10L^3-10L^2x+5Lx^2-x^3)
    \end{equation}
\end{enumerate}
In all cases the Young’s modulus is $E=1$ and the Poisson ratio is $\nu=0$, and the result is shown in \cref{fig:2D_Bending}. 
A good match can be observed.

\subsubsection{2D/3D stress concentration}
\begin{figure}[tbhp]
    \centering
    \includegraphics[width=\textwidth]{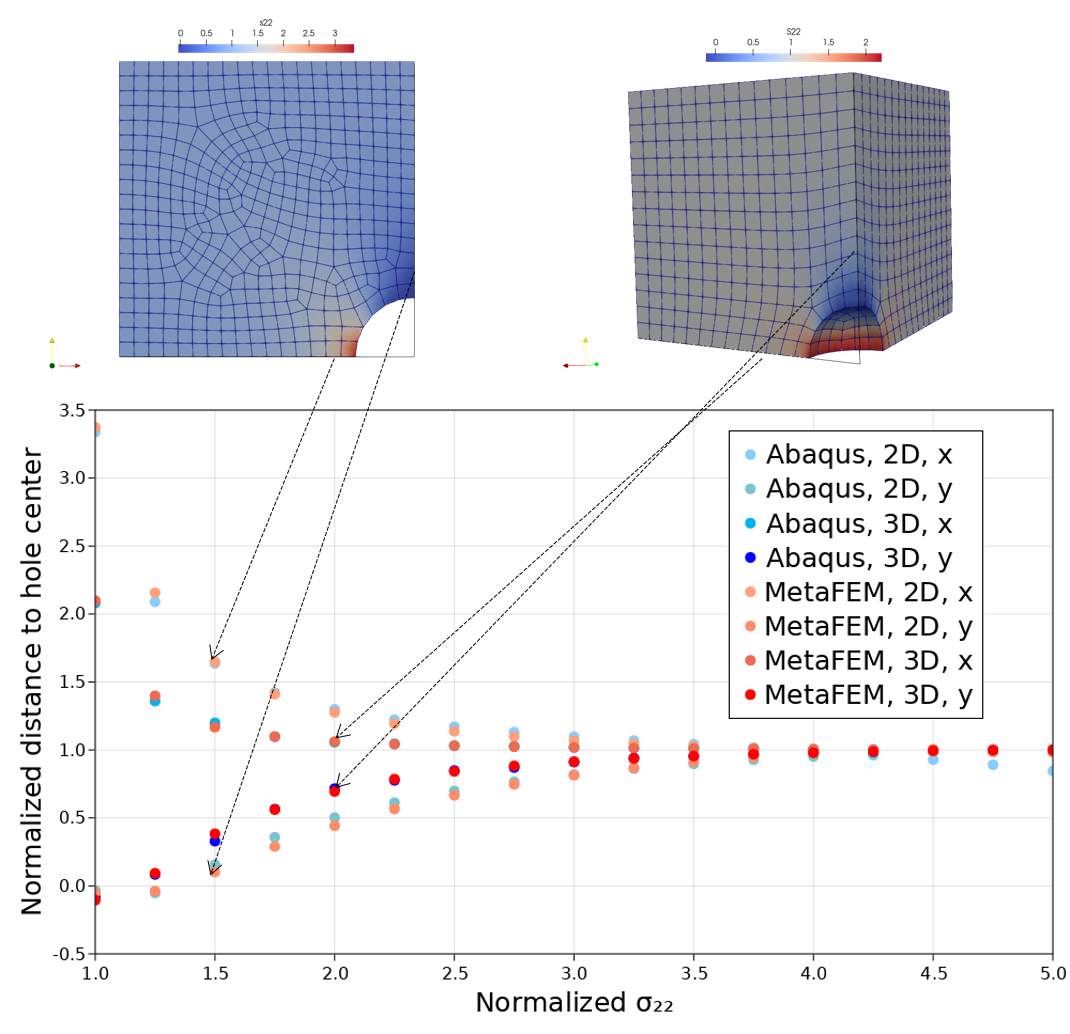}
    \caption{2D/3D stress concentration. Top left: Mesh and $\sigma_{22}$ distribution for 2D; Top right: Mesh and $\sigma_{22}$ distribution for 3D;
    Bottom: Comparison of $\sigma_{22}$ along two axes.}
    \label{fig:3D_SC}
\end{figure}

A square/cube of side length $L=10\,m$ with a circular/spherical hole of radius $r=1\,m$ under uniaxial tension $\sigma_0$ along the y axis is 
simulated by both MetaFEM 
and Abaqus with 440/3,375 quadratic serendipity quadrilateral/brick elements and 1,399/15,645 mesh nodes, in total $2\times2=4$ simulations. 
The material is chosen as the general steel, $E = 210\times10^9\,Pa, \nu = 0.3$.
Only $\frac{1}{4}$/$\frac{1}{8}$ of the sample is simulated 
to take advantage of symmetry.
$\sigma_{22}$ distribution along the $x$ and $y$ axes are compared in \cref{fig:3D_SC}, where we recover the well-known 2D/3D stress concentration factors 3 and 2 on the 
intersection of the circular/spherical surface and the x axis, while the distributions calculated by MetaFEM match well with the Abaqus results. 

\subsection{Incompressible flow}
The mathematical formulation of the steady incompressible Navier-Stokes (NS) equation in FEM
with the Streamline Upwind/Petrov-Galerkin (SUPG) stabilization is given by:
\begin{equation}
    \nonumber \text{variable}\quad u_i,p,\qquad\text{parameters}\quad\rho,\mu,u^W_i,\tau^b,\tau^c,\tau^m
\end{equation}
\begin{equation}
    \nonumber Rc\coloneqq u_{k,k},\qquad {Rm}_i\coloneqq\rho u_ku_{i,k}+p_{,i}-\mu u_{i,kk}
\end{equation}
\begin{align}
    \nonumber \overbrace{-\rho(u_{i,j},u_iu_j)-(u_{i,i},p)+(p,u_{i,i})+\mu(u_{i,j},u_{i,j})}^{NS}+\\
    \overbrace{\tau^m\rho(u_{i,j},{Rm}_iu_j)+\tau^m(p_{,i},{Rm}_i)+\tau^c(u_{i,i},Rc)}^{SUPG} \\
    \nonumber =0,\qquad in\quad\Omega
\end{align}

\begin{alignat}{4}
    \nonumber (u_i,pn_i)-\mu(u_i,u_{i,j}n_j)={}&0,\qquad on\quad\partial\Omega\\
    \nonumber \rho(u_i,u^w_iu^w_jn_j)+(p,(u^w_i-u_i)n_i)\qquad&\\
    \nonumber +\mu(u_{i,j},(u^w_i-u_i)n_j)+\tau^b\rho(u_i,u_i-u^w_i)={}&0,\qquad on\quad(\partial\Omega)_{inflow}\\
    \nonumber \rho(u_i,u_iu_jn_j)={}&0,\qquad on\quad(\partial\Omega)_{outflow}\\
    (p,-u_in_i)+\mu(u_{i,j},-u_in_j)+\tau^b\rho(u_i,u_i)={}&0,\qquad on\quad(\partial\Omega)_{fix}
\end{alignat}
where $u_i,p,\rho,\mu$ are the velocity, pressure, density and the dynamic viscosity respectively.
The code for physics is given by:
\begin{jllisting}
@Def begin
    Rc = u{k;k}
    Rm{i} = ρ * u{k} * u{i;k} + p{;i} - μ * u{i;k,k} 
end

@Def begin
   NS_domain_BASE = - ρ * Bilinear(u{i;j}, u{i} * u{j}) - Bilinear(u{i;i}, p) + 
   Bilinear(p, u{i;i}) + μ * Bilinear(u{i;j}, u{i;j})

   NS_domain_SUPG = τᵐ * ρ * Bilinear(u{i;j}, Rm{i} * u{j}) + 
   τᵐ * Bilinear(p{;i}, Rm{i}) + τᶜ * Bilinear(u{i;i}, Rc)

   NS_boundary_BASE = Bilinear(u{i}, p * n{i}) - μ * Bilinear(u{i}, u{i;j} * n{j})

   NS_boundary_INFLOW = ρ * Bilinear(u{i}, uʷ{i} * uʷ{j} * n{j}) + 
   Bilinear(p, (uʷ{i} - u{i}) * n{i}) + μ * Bilinear(u{i;j}, (uʷ{i} - u{i}) * n{j}) + 
   τᵇ * ρ * Bilinear(u{i}, u{i} - uʷ{i})

   NS_boundary_OUTFLOW = ρ * Bilinear(u{i}, u{i} * u{j} * n{j})  

   NS_boundary_FIX = Bilinear(p, - u{i} * n{i}) + μ * Bilinear(u{i;j}, - u{i} * n{j}) + 
   τᵇ * ρ * Bilinear(u{i}, u{i})
end

@Def begin
   NS_domain = NS_domain_BASE + NS_domain_SUPG
   NS_boundary_inflow = NS_boundary_BASE + NS_boundary_INFLOW
   NS_boundary_outflow = NS_boundary_BASE + NS_boundary_OUTFLOW
   NS_boundary_fix = NS_boundary_BASE + NS_boundary_FIX
end
\end{jllisting}

\subsubsection{2D lid-driven cavity}
\begin{figure}[tbhp]
    \centering
    \includegraphics[width=\textwidth]{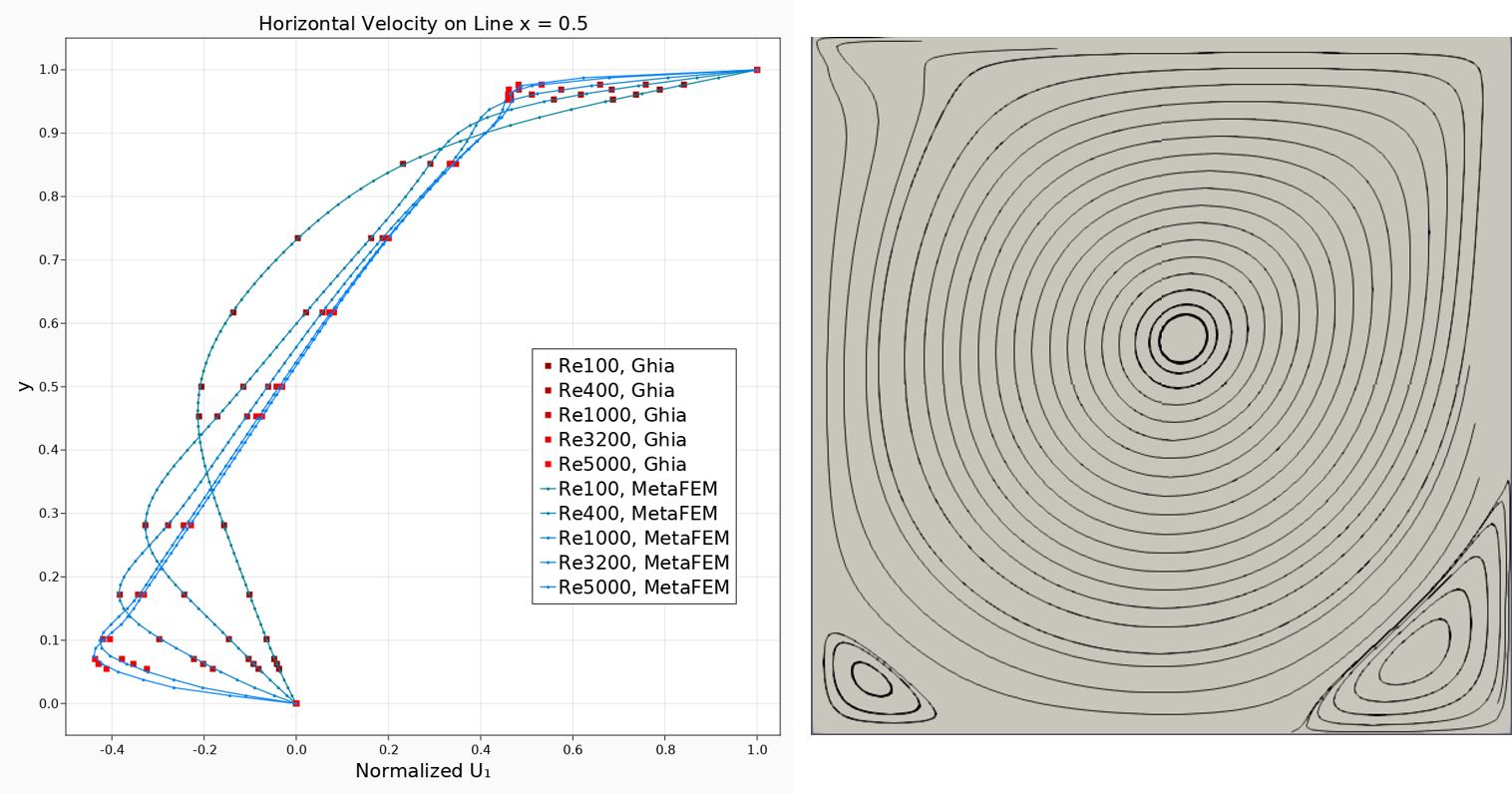}
    \caption{The lid-driven cavity flow. Left: Horizontal velocity comparison along the middle vertical line $x = 0.5$; 
    Right: Streamline plot for $Re=$ 1,000}
    \label{fig:2D_CavityFlow}
\end{figure}

A cavity of side length $L=1$ full of fluid $\rho=1,\mu=1$ has its left, bottom and right wall fixed and the top wall moving rightward with a velocity 
$u^w$ determined by the Reynold number, $u^w=Re\frac{\mu}{L\rho}$, 
is simulated for cases with $Re=$ 100, 400, 1,000, 3,200 and 5,000 with $40\times40$
quadratic serendipity square elements (dx = 0.025). 
The horizontal velocity on the middle vertical line $x = 0.5$ is compared between MetaFEM and the result of Ghia \cite{cavity_1982}, 
as shown in \cref{fig:2D_CavityFlow}. The Paraview streamline plot for the example case $Re=$1,000 is also provided.

\subsection{3D cylinder flow}
A cuboid channel of length $L=2.5$, width and depth $h=0.41$ with a fixed z-axis oriented cylinder obstacle of radius $r=0.05$ 
is positioned with its center at $(x,y) = (0.2,0.2)$. The 4 channel sides ($y=0,z=0,y=0.41,z=0.41$) are fixed while a steady flow 
($\rho=$ 1,000, $\mu=$ 1) with a fully developed parabolic inflow profile at the left entry $x=0$ with maximum velocity $U=0.45$ ($Re=20$) 
flows out to right outlet $x=2.5$. The detailed inflow profile is:
\begin{equation}
    u^w(0,y,z)=(\frac{16U(h-y)(h-z)yz}{h^4},0,0)
\end{equation}

\begin{figure}[tbhp]
    \centering
    \includegraphics[width=\textwidth]{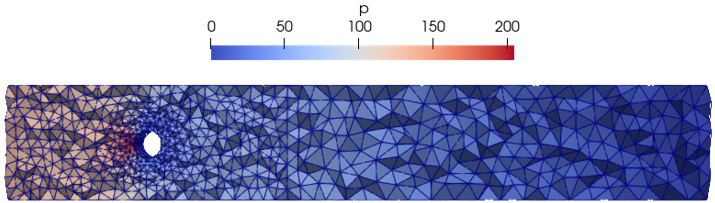}
    \caption{The mesh for the cylinder flow, colored with pressure distribution}
    \label{fig:3D_CF_Pclip}
\end{figure}

The case is simulated in both COMSOL and MetaFEM with 28,468 quadratic simplex elements and with 41,202 mesh nodes. The mesh is shown 
\cref{fig:3D_CF_Pclip}.

\begin{figure}[tbhp]
    \begin{subfigure}{\textwidth}
        \centering
        \includegraphics[width=1.\linewidth]{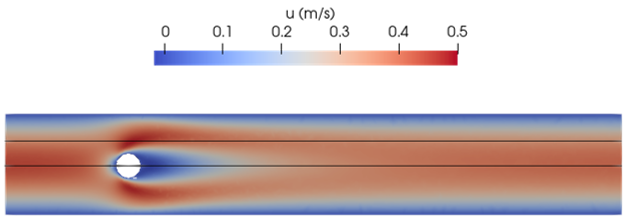}   
        \caption{The plan $z=0.2$ for the sample lines, colored with horizontal velocity distribution}
    \end{subfigure}

    \begin{subfigure}{0.9\textwidth} 
        \centering
        \includegraphics[width=0.9\linewidth]{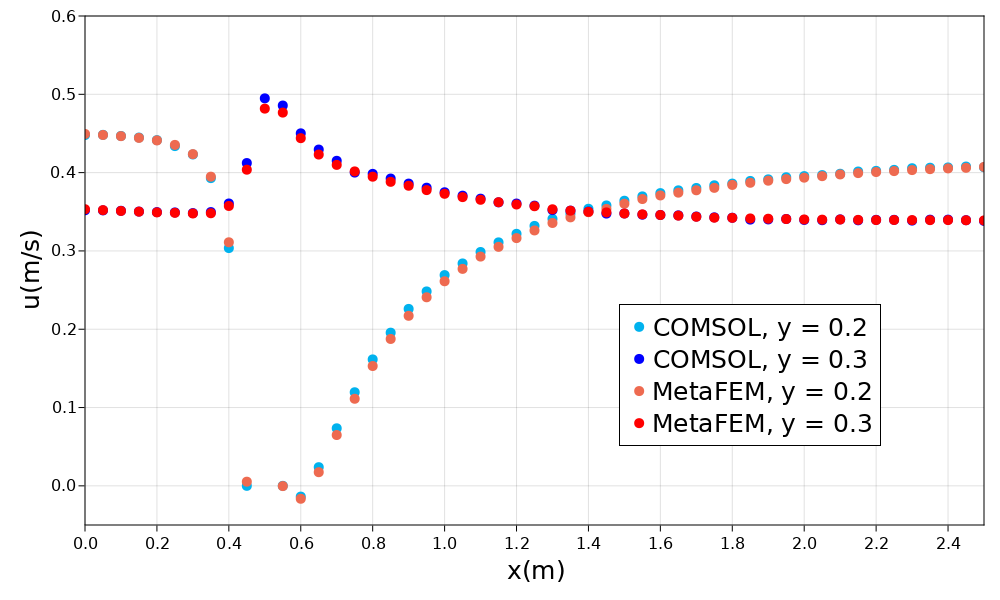} 
        \caption{Horizontal velocity distribution comparison.}
        \label{fig:CF_u}
    \end{subfigure}

    \begin{subfigure}{0.9\textwidth}
        \centering
        \includegraphics[width=0.9\linewidth]{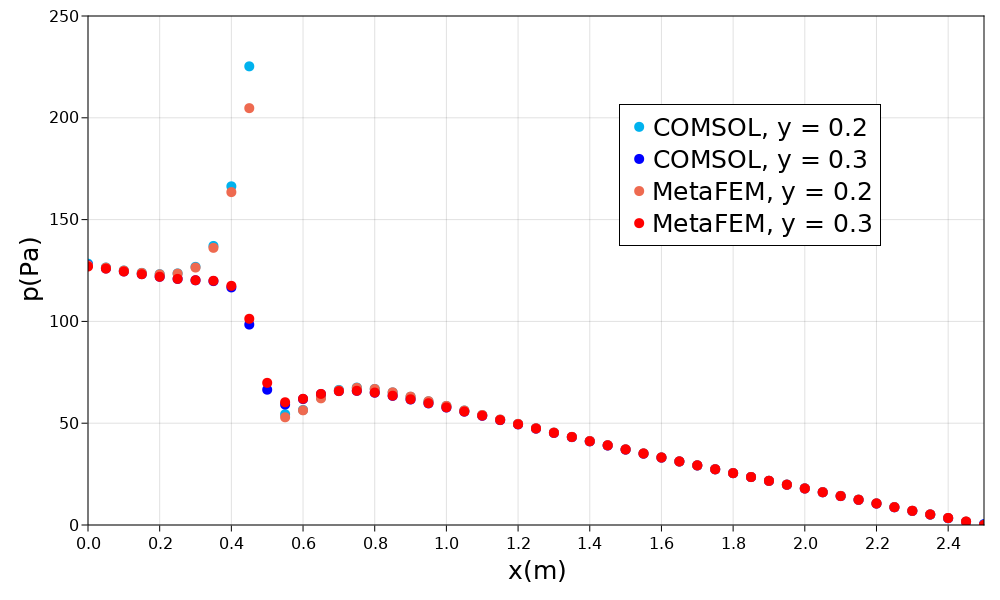}  
        \caption{Pressure comparison.}
    \end{subfigure}
    \caption{The 3D cylinder flow}
    \label{fig:3D_CF_Rest}
\end{figure}

For further verification, we compare the horizontal velocity $u_1$ and pressure $p$ along two horizontal lines 
$y=0.2, z=0.2$ and $y=0.3, z=0.2$, as in \cref{fig:3D_CF_Rest}, where the good match indicates that MetaFEM is able to correctly derive and assemble the physics 
of incompressible laminar flow.

\section{Discussions}

\label{sec:conclusions}
This paper proposes MetaFEM, a generic finite element
solver with original formulations of the theory, algorithm and of the final usable software.
Its theoretical/academic contribution comes from the meta-expression formulation and the rewriting system design.
Its practical/engineering contribution comes from a compact skeleton software which is able to simulate various PDE systems like thermal conduction in solids,
linear elasticity and incompressible flow while depending only on basic/generic Julia packages. 

The initial motivation for MetaFEM, which is still one of the final goals, is to provide a fast evaluation for practical manufacturing processes, 
e.g., metal forming, casting or additive manufacturing, where both the domain physics and the boundary conditions can be highly customized.
However, despite of the fact that most cases can be described by the meta-expression formulation, there are still two major limitations 
preventing MetaFEM from describing practical processes: (1) FEM is memory-intense, 
GPU-accelerated FEM is apparently even more memory-intense. 
A generic distributed design is certainly needed but has not been implemented.
(2) Practical processes usually lead to large deformations, where traditional static mesh systems will suffer from mesh distortion.
Under such circumstances, neither the Arbitrary Lagrangian-Eulerian (ALE) formulation or remeshing rules can fully alleviate the problem.
A cutcell mesh system may be a fundamental solution but has not been implemented either. 
We left the above two points as major future works.

\end{document}